\documentclass[11pt,reqno]{amsart}
\usepackage{amssymb,amsthm,amsfonts,amsmath,amscd,amssymb,epsfig,psfrag}
\usepackage{calrsfs}
\usepackage[all]{xy}

\usepackage{amssymb,bbm,graphicx,epsfig,psfrag,epic}
\usepackage{enumerate,array,hyperref,graphicx,color,amsthm,bbm,mathrsfs,eepic,latexsym}
\usepackage{mathrsfs}
\theoremstyle{plain}
\newtheorem{theorem}{Theorem}[section]
\newtheorem{lemma}[theorem]{Lemma}

\newtheorem*{remark*}{Remark}
\newtheorem*{remarks*}{Remarks}

\newtheorem{example}[theorem]{Example}

\newtheorem*{example*}{Example}                           
\newtheorem*{examples*}{Examples}
\newtheorem*{examplecon*}{Example (continued)}

\newtheorem{question}[theorem]{Question}

\DeclareMathAlphabet{\pazocal}{OMS}{zplm}{m}{n}

\newcommand{\NN}{{\mathbb{N}}}
\newcommand{\ZZ}{{\mathbb{Z}}}
\newcommand{\RR}{{\mathbb{R}}}
\newcommand{\CC}{{\mathbb{C}}}
\newcommand{\QQ}{{\mathbb{Q}}}
\newcommand{\FF}{{\mathbb{F}}}

\newcommand{\m}{{\bf m}}

\newcommand{\Crit}{{\rm Crit}\,}

\newcommand{\loc}{{\rm loc}}
\newcommand{\ev}{{\rm ev}}

\newcommand{\HF}{{\rm HF}}

\def\mcA{{\mathcal{A}}}
\def\mcP{{\mathcal{P}}}

\newcommand{\proofend}{\hspace*{\fill} $\Box$\\}

\newcommand{\pr}{{\mathrm pr}}

\def\1{\:\!}
\def\2{\;\!}
\def\ni{\noindent}

\def\m{\medskip}
\def\s{\smallskip}
\def\pp{\partial}
\def\id{\operatorname{id}}
\def\pr{\operatorname{pr}}

\def\cuplength{\operatorname{cuplength}}
\def\tot{\operatorname{tot}}
\def\CP{\operatorname{\mathbbm{C}P}}

\begin{document}

\title[Iterated graph construction and Hamiltonian delay equations]{An iterated graph construction and periodic orbits of Hamiltonian delay equations}

\author{Peter Albers}
\address{Peter Albers, Mathematisches Institut, Universit\"at Heidelberg}
\email{palbers@mathi.uni-heidelberg.de}

\author{Urs Frauenfelder}
\address{Urs Frauenfelder, Mathematisches Institut, Universit\"at Augsburg}
\email{urs.frauenfelder@math.uni-augsburg.de}

\author{Felix Schlenk}  
\address{Felix Schlenk, Institut de Math\'ematiques, Universit\'e de Neuch\^atel}
\email{schlenk@unine.ch}

\keywords{Delay equations, Hamiltonian systems, Floer homology}

\date{\today}
\thanks{2000 {\it Mathematics Subject Classification.}
Primary 34K13, Secondary~58E05, 58F05}

\begin{abstract}
According to the Arnold conjectures and Floer's proofs, 
there are non-trivial lower bounds for the number of periodic solutions
of Hamiltonian differential equations on a closed symplectic manifold whose symplectic form vanishes on spheres.
We use an iterated graph construction and Lagrangian Floer homology to 
show that these lower bounds also hold for certain Hamiltonian delay equations.
\end{abstract}

\maketitle
\tableofcontents

\section{Introduction and main results} 
\label{s:intro}

\ni
{\bf Delay equations.}
An ordinary differential equation (ODE) on a manifold~$M$
is a problem of the form 
$$
\dot x (t) \,=\, X_t (x(t)) 
$$
where~$X_t$ is a time-dependent vector field on~$M$.
A delay differential equation (DDE) on~$M$ is a problem of the form 
\begin{equation} \label{e:delay}
\dot x (t) \,=\, \mathcal{X} (x_{\leq t}) 
\end{equation}
where $x_{\leq t} = \{ x(s) \colon s \leq t \}$ is the ``past before time~$t$'' of the curve~$x$
and where $\mathcal{X}$ associates with each curve $x_{\leq t}$, $t \in \RR$,
a vector in~$T_{x(t)}M$.
In a DDE, the velocity vector $\dot x (t)$ therefore does not only depend on the
instantaneous position~$x(t)$, but also on positions of~$x$ in the past.
Note that if $x$ is a periodic orbit,  
then $x_{\leq t}$ is the whole orbit~$x$. 
Delay equations form a vast topic;
we refer to~\cite{HaVL93} for the general theory
and to~\cite{Erneux09} for a wealth of DDEs arising in applications.

The simplest delay equations on~$\RR^d$ are of the form
\begin{equation} \label{e:delayX}
\dot x(t) \,=\, \sum_{j=1}^N X_t^j \bigl( x(t-\tau_j) \bigr)
\end{equation}
where $X_t^j$ are time-dependent vector fields on~$\RR^d$ and $0 \leq \tau_1 < \dots < \tau_N$ 
are the delay times.
On a general manifold~$M$, one has to ensure that both sides of~\eqref{e:delayX} live in the same tangent space, 
that is, the vector fields~$X_t^j$ are of the form 
\begin{equation} \label{e:delayY1}
\dot x(t) \,=\, \sum_{j=1}^N f_t^j \bigl( x(t-\tau_j) \bigr) \, Y_t^j  \bigl( x(t) \bigr)
\end{equation}
for time-dependent functions $f_t^j$ and vector fields $Y_t^j$ on~$M$.

Assume now that $(M,\omega)$ is a closed symplectic manifold such that the cohomology class $[\omega]$ vanishes on spherical homology classes: $[\omega]|_{\pi_2(M)} =0$.
%
%
Also assume that $X_H$ is the Hamiltonian vector field of a function $H \colon M \times S^1 \to \RR$,
where $S^1 = \RR / \ZZ$.
Then by the proofs of Arnold's celebrated conjectures
(see for instance~\cite[Chapter~11]{mcduff-salamon17}) 
the number of 1-periodic solutions of the Hamiltonian ODE $\dot x(t) = X_H(x(t))$ is at least 
$\cuplength (M)+1$, and at least $\dim H_* (M;\ZZ_2)$ generically 
(namely if the graph of $\phi_H^1$ in $M \times M$ intersects the diagonal transversally).
Here, $\cuplength (M)$ is the maximal length of a non-vanishing product of elements of positive degree in the cohomology ring~$H^*(M;\ZZ_2)$,
and $\dim H^* (M;\ZZ_2)$ is the sum of the Betti numbers of~$M$ with respect to $\ZZ_2$ coefficients.  

\begin{question}
Do the Arnold conjectures generalize to delay equations?
\end{question}

To make the question meaningful, one should first decide what a Hamiltonian delay equation should be. 
In~\cite{AFS17.3}, we give one possible answer to the sub-question what a periodic solution of a Hamiltonian 
delay equation should be, by characterising such solutions variationally. 
In this note we prove the Arnold conjectures for a special class of Hamiltonian delay equations by
using classical tools from the theory of $J$-holomorphic curves and Floer homology, 
and by an iterated graph construction. 
The delay equations we can deal with here in particular have only a finite number of delay times.
We believe that the Arnold conjectures hold true for a far more general class of Hamiltonian delay
equations, cf.\ \cite{AFS17.3}. 
For this one should construct a non-local Floer homology, in which the
gradient flow lines will not anymore be solutions of an elliptic PDE on a finite dimensional symplectic manifold, 
but an ODE on an infinite dimensional scale manifold.  
For a detailed description of this kind of equations we refer to~\cite{AFS17.1},
where a first step of this construction (namely compactness of the space of solutions of non-local Floer equations) 
was taken.
The results of this paper, that rely only on classical Floer theory, 
indicate that such a more general Floer theory, that would imply Theorem~\ref{t:main} below as a special case, 
does exist.

\m \ni
{\bf Delay equations by means of an iterated graph construction.} 
In accordance with A.~Weinstein's dictum  from~\cite{Wei81} that
``everything is a Lagrangian submanifold'',
many results on periodic orbits of Hamiltonian systems on a symplectic manifold~$(M,\omega)$
can be deduced from results on Lagrangian intersections. 
Indeed, for every Hamiltonian function $H \colon M \times S^1 \to \RR$
the intersections of the diagonal~$\Delta$ in $(M \times M, \omega \oplus - \omega)$ with the graph
$\{ (\phi_H^1(z), z)\}$ of the Hamiltonian diffeomorphism~$\phi_H^1$ 
are in bijection with the 1-periodic orbits of the Hamiltonian flow of~$H$.

To give this correspondence a more precise dynamical interpretation, 
we consider the loop space $\mcP_0 = W^{1,2}(S^1,M)$
and the path space $\mcP_1 = \left\{ v \in W^{1,2} \left( [0,1], M \times M \right) \mid v(0), v(1) \in \Delta \right\}$.
The map 
$$
\psi \colon \mcP_0 \to \mcP_1, \quad  \psi (v)(t) = (v(t), v(0))
$$
is an embedding.
Given $H \colon M \times S^1 \to \RR$, define 
$\widetilde H \colon M \times M \times S^1 \to \RR$ by
$$
\widetilde H_t (z_1,z_2) \,=\, H_t(z_1).
$$
Then the flow-lines of $\widetilde H$ are of the form $(\phi_H^t(z_1),z_2)$, and so the
Hamiltonian chords of $\widetilde H$ that start and end at time~1 on the diagonal 
correspond to $1$-periodic orbits of~$H$ on~$M$.

One may now wonder what more general Hamiltonian functions on $M \times M$ can tell us about dynamical systems on~$M$.
The flow lines of such more general Hamiltonians may, however, not be in the image of~$\psi$,
and thus cannot be pulled back to~$M$ by~$\psi$.
Following~\cite[\S 5.2]{BPS03} we therefore consider another map
$\Psi \colon \mcP_0 \to \mcP_1$, given by
\begin{equation*} 
\Psi(v)(t) = \Big( v\big(\tfrac{t}{2}\big),v\big(1-\tfrac{t}{2}\big) \Big), \quad t \in [0,1] ,
\end{equation*}
which in contrast to $\psi$ is a diffeomorphism, with inverse
\begin{equation*} 
\Psi^{-1}(w)(t) = \left\{
\begin{array}{ll}
w_1(2t),   & t \in \big[ 0,\tfrac{1}{2} \big], \\ [0.3em]
w_2(2-2t), & t \in \big[ \tfrac{1}{2},1 \big].
\end{array}\right.
\end{equation*}
Under the ``change of variables'' $\Psi$,
the Hamiltonian loops of $H \colon M \times S^1 \to \RR$ correspond to Hamiltonian chords of
$$
\widetilde H \colon M \times M \times S^1 \to \RR, \quad
\widetilde H_t(z_1,z_2) = \tfrac 12 H_{\frac t2}(z_1) + \tfrac 12 H_{1-\frac t2}(z_2) .
$$
For other functions $K \colon M \times M \times S^1 \to \RR$,
however, Hamiltonian chords are pulled back to periodic orbits on~$M$ that solve
a {\it delay equation}.
For instance, taking $K$ in product form, $K_t(z_1,z_2) = \tfrac 12 \1 F_t(z_1) \, G_t(z_2)$,
Hamiltonian chords of $K$ are pulled back to $1$-periodic orbits~$v$ solving
\begin{equation*} 
\dot v (t) \,=\, 
\left\{
\begin{array}{ll}
G_{2t} (v(1-t)) \,X_{F_{2t}} (v(t)),       & t \in \bigl[ 0,\tfrac{1}{2} \bigr], \\ [0.3em]
F_{2-2t} (v(1-t)) \, X_{G_{2-2t}} (v(t)) , & t \in \bigl[ \tfrac{1}{2},1 \bigr].
\end{array}\right.
\end{equation*}
While the term $v(1-t)$ may look a bit awkward at a first glance, since time is running backwards, 
along 1-periodic orbits the above equation is an honest delay equation of the form~\eqref{e:delay}, 
because along such orbits we have $x(t+k) = x(t)$ for all~$k \in \ZZ$ and so there is no difference 
between the future and the past.

To obtain periodic solutions of more general Hamiltonian delay equations on~$M$, 
and in particular such of the classical form~\eqref{e:delayY1}, we iterate the above graph construction:
Set $(M_1, \omega_1) = (M \times M, \omega \oplus -\omega)$ and
define the symplectic manifold
$$
(M_2,\omega_2) := (M_1 \times M_1, \omega_1 \oplus -\omega_1)
 = (M^4, \omega \oplus -\omega \oplus -\omega \oplus \omega), 
$$
and in $(M_2, \omega_2)$ take two Lagrangian submanifolds, 
namely the product of the diagonal in~$M_1$ with itself:
$$
\Delta_2^0 := \Delta \times \Delta \,\subset\, M_2,
$$
and the diagonal
$$
\Delta_2^1 := \{(z,z) : z \in M_1\} \,\subset\, M_2.
$$
Then the same formula as for $\Psi$ above yields a diffeomorphism
$\Psi_1$ from~$\mathcal{P}_1$ to the space
$$
\mathcal{P}_2 = \big\{ v \in W^{1,2}([0,1],M_2) : 
            v(0) \in \Delta^0_2,\,\,v(1) \in \Delta^1_2 \big\} .
$$
Pulling back Hamiltonian chords of, say, 
$$
K_t(z) \,=\, F_t^1(z_1) \, F_t^4(z_4) + F_t^2(z_2) \, F_t^3(z_3) 
$$
on $M_2$ by $\Psi^1 = \Psi_1 \circ \Psi$ we then get 1-periodic orbits on~$M$
that solve a delay equation with delay time~$\frac 12$.

Iterating the graph construction, we obtain for every $n \geq 1$ 
a diffeomorphism $\Psi^n$ from $\mcP_0$ to the space~$\mcP_n$ of paths in~$M_n := M^{2^n}$
that end on the diagonal~$\Delta_n^1$ of $M_n = M_{n-1} \times M_{n-1}$
and start on $\Delta_n^0 := \Delta_{n-1}^0 \times \Delta_{n-1}^1$.
From this diffeomorphism we get 1-periodic orbits on~$M$ that solve delay equations
with time delays $\frac{1}{2^n}, \dots, \frac{2^n-1}{2^n}$,
and pre-composing~$\Psi$, $\Psi_1, \dots$ with a diffeomorphism of the circle or the interval,
we get delay equations with various other delay times,
see Sections~\ref{s:delay} and~\ref{s:general}.

The existence of many Hamiltonian chords on $\mcP_n$ is guaranteed by 
a recent cuplength estimate of Albers--Hein~\cite{AlHe16} 
and by Pozniak's theorem on clean Lagrangian intersections~\cite{Po99}.
From this we obtain

\begin{theorem} \label{t:main}
Let $(M,\omega)$ be a closed symplectic manifold with $[\omega] |_{\pi_2(M)} =0$,
and fix a function $K \colon M_n \times S^1 \to \RR$.
Let $\Psi^n \colon \mcP_0 \to \mcP_n$ be the diffeomorphism obtained by the iterated graph construction, 
or any of its variants induced by pre-composing with diffeomorphisms of the circle or the interval.
Then the number of contractible 1-periodic orbits on~$M$ that solve the Hamiltonian delay equation 
induced by~$K$ via pullback under~$\Psi^n$ is
\begin{itemize}
\item[(i)] at least $\cuplength (M) +1$;

\s
\item[(ii)]
at least $\dim H_* (M;\ZZ_2)$ 
if $\phi_K^1 (\Delta_n^0)$ intersects~$\Delta_n^1$ transversally.
\end{itemize}
\end{theorem}

\ni 
Note  that the transversality assumption on $K$ in (ii) is $C^\infty$-generic in the space of functions $M_n \times S^1 \to \RR$.

\m \ni
{\bf Comparison with previous results.}
We conclude this introduction with comparing our method and result with previous ones.

The search for periodic solutions of (non-linear) DDEs, 
that has been an important topic in the study of DDEs since the~1960s, 
is much harder than for ODEs:
Most results are for very special classes of autonomous DDEs on~$\RR$, 
and just one or two periodic orbits are found.
Existence results for periodic orbits of DDEs on~$\RR^d$, or even on manifolds, 
are even scarcer, see \cite{BCF09},
~\cite[\S 11]{HaVL93},~\cite{Wal14} and the references therein.
The methods used are bifurcation and fixed point theorems, and index methods. 
In addition to these tools,
there are a few tricks that for some DDEs convert the problem of finding periodic solutions to the problem of finding certain solutions of related~ODEs. 
One is the ``chain trickery'' going back to~\cite{Far73}, 
see also~\cite{BoHa00}.
Another one is the Kaplan--Yorke method from~\cite{KaYo74}.
Our graph construction adds one more trick to this lists.

While the search for periodic orbits of {\it Hamiltonian} ODEs is a topic 
with a rich tradition and many profound results, 
much less is known on the existence of periodic orbits of Hamiltonian DDEs.
A Hamiltonian DDE on~$\RR^{2n}$ is, for instance, 
of the form~\eqref{e:delayX} with $X_t^j = J \nabla H_t^j$, 
where $J$ is the usual complex structure on $\RR^{2n} \cong \CC^n$.
For the special case
$$
\dot x (t) \,=\, J 
\bigl( \nabla G_t(x(t-\tau)) + \nabla G_t(x(t-2\tau)) + \dots + \nabla G_t(x(t-N\tau)) \bigr)
$$
with $G_t$ of periodic~$\tau$ and meeting suitable growth conditions, 
Liu~\cite{Liu12} proved the existence of one $N \tau$-periodic solution, 
and, under a generic assumption on this solution, 
of two or even three $N \tau$-periodic solutions. 
The main idea of his proof, that can be traced back to the Kaplan--Yorke method, 
was one inspiration for our graph trick.

In summary, 
both our methods and results are different to previous ones: 
We use the {\it graph trick}\, and {\it Floer homology}\, to prove a 
{\it multiplicity}\, result for periodic solutions of certain {\it Hamiltonian}\, 
delay equations on {\it manifolds}.

\m
The paper is organized as follows. 
In the next section we describe the iterated graph construction in detail.
In Section~\ref{s:action} we investigate how the action functionals of classical mechanics on~$\mcP_0$
and~$\mcP_n$ relate unter~$\Psi^n$.
In Sections~\ref{s:delay} and~\ref{s:general} we give many examples of Hamiltonian delay equations
for periodic orbits that are obtained by pulling back Hamiltonian chords on $\mcP_n$ by~$\Psi^n$
and its variants.
In Section~\ref{s:proof} we prove Theorem~\ref{t:main},
and in the last section we discuss improvements of Theorem~\ref{t:main} in various directions.

\subsubsection*{Acknowledgment}
FS cordially thanks Augsburg University for its warm hospitality in the autumn of~2017.
We are grateful to Felix Schm\"aschke for his explanations on orienting the space of Floer strips.
PA is supported by DFG CRC/TRR 191, UF is supported by DFG FR/2637/2-1, and FS supported by SNF grant 200020-144432/1


\section{The iterated graph construction} \label{s:graph}
Assume that $(M,\omega)$ is a symplectic manifold. 
Abbreviate by $S^1 = \RR/\ZZ$ the circle and let  
$$
\mathcal{P}_0 = W^{1,2}(S^1,M)
$$
be the space of free loops on $M$ of Sobolev class~$W^{1,2}$. 
The product manifold 
$$
(M_1, \omega_1) := 
(M \times M, \omega \oplus -\omega)
$$ 
is again a symplectic manifold, and the diagonal
$$
\Delta_1 = \big\{ (z,z) : z \in M \big\} \subset M_1
$$
is a Lagrangian submanifold canonically diffeomorphic to~$M$ 
via the map $(z,z) \mapsto z$. 
Abbreviate by
$$
\mathcal{P}_1 = \big\{ v \in W^{1,2}([0,1],M_1) : 
               v(0), v(1) \in \Delta_1 \big\}
$$
the space of paths of class $W^{1,2}$ in $M \times M$ that start and end 
at the diagonal. 
Note that the boundary condition makes sense since the elements of $W^{1,2}([0,1],M_1)$
are continuous.

The Hilbert manifolds~$\mathcal{P}_0$ and~$\mathcal{P}_1$ are diffeomorphic 
via the diffeomorphism
$\Psi \colon \mathcal{P}_0 \to \mathcal{P}_1$ given by
\begin{equation} \label{e:Psi}
\Psi(v)(t) = \Big( v\big(\tfrac{t}{2}\big),v\big(1-\tfrac{t}{2}\big) \Big), \quad t \in [0,1] .
\end{equation}
Indeed, since $v \in \mathcal{P}_0$ is a continuous loop it holds that $v(0)=v(1)$, so that
$$
\Psi(v)(0) = \big( v(0), v(1) \big) = \big( v(0),v(0) \big) \in \Delta_1 .
$$
Moreover,
$\Psi(v)(1) = \big( v\big( \tfrac{1}{2} \big), v \big( \tfrac{1}{2} \big) \big) \in \Delta_1$,
whence $\Psi(v)$ satisfies the required boundary conditions, i.e., 
$\Psi(v) \in \mathcal{P}_1$. 
The inverse $\Phi \colon \mathcal{P}_1 \to \mathcal{P}_0$ of~$\Psi$ can be explicitly 
written down:
For $w=(w_1,w_2) \in \mathcal{P}_1$ it is given by
\begin{equation} \label{e:Phi}
\Phi(w)(t) = \left\{
\begin{array}{ll}
w_1(2t),   & t \in \big[ 0,\tfrac{1}{2} \big], \\ [0.3em]
w_2(2-2t), & t \in \big[ \tfrac{1}{2},1 \big].
\end{array}\right.
\end{equation}
Since $w \in \mathcal{P}_1$ it holds that $w_1(0)=w_2(0)$ and $w_1(1)=w_2(1)$, 
which guarantees that $\Phi(w)(0)=\Phi(w)(1)$ and that $\Phi(w)$ is continuous at $t=\tfrac{1}{2}$,
so that $\Phi(w)$ belongs to $\mathcal{P}_0$.
Further, $\Phi \circ \Psi = \id_{\mathcal{P}_0}$ and $\Psi \circ \Phi = \id_{\mathcal{P}_1}$.
We note that it is important that we work with loops and paths of class $W^{1,2}$.
For spaces $W^{k,2}$ of higher regularity, $\Psi$ would not be surjective
and $\Phi$ would not take values in the space of $W^{k,2}$-loops.

The next iteration step proceeds as follows. Define the symplectic manifold
$$
(M_2,\omega_2) := (M_1 \times M_1, \omega_1 \oplus -\omega_1)
 = (M^4, \omega \oplus -\omega \oplus -\omega \oplus \omega).
$$
In $(M_2, \omega_2)$ we introduce two Lagrangian submanifolds, 
namely the product of the diagonal in~$M_1$ with itself:
$$
\Delta_2^0 := \Delta_1 \times \Delta_1 \,\subset\, M_2,
$$
and the diagonal
$$
\Delta_2^1 := \{(z,z) : z \in M_1\} \,\subset\, M_2.
$$
Note that both Lagrangians $\Delta_2^0$ and $\Delta_2^1$ are canonically 
diffeomorphic to $M_1 = M \times M$. Their intersection 
$$
\Delta^0_2 \cap \Delta^1_2 = \{(z,z,z,z) : z \in M\} \subset M_2 = M^4
$$
is the total diagonal in the fourfold product~$M^4$ 
and is therefore canonically diffeomorphic to~$M$.
We introduce 
$$
\mathcal{P}_2 = \big\{ v \in W^{1,2}([0,1],M_2) : 
            v(0) \in \Delta^0_2,\,\,v(1) \in \Delta^1_2 \big\} ,
$$
the space of paths in $M_2$ of Sobolev class~$W^{1,2}$ which start 
at~$\Delta^0_2$ and end in~$\Delta^1_2$. 
This Hilbert manifold is diffeomorphic to~$\mathcal{P}_1$ 
and therefore also to~$\mathcal{P}_0$. 
A diffeomorphism $\Psi_1 \colon \mathcal{P}_1 \to \mathcal{P}_2$ 
is given by the same formula as for~$\Psi$, 
namely
\begin{equation*} 
\Psi_1(v)(t) = \Big( v\big( \tfrac{t}{2} \big), v\big( 1-\tfrac{t}{2} \big) \Big), \quad t \in [0,1] .
\end{equation*}
The inverse $\Phi_1 \colon \mathcal{P}_2 \to \mathcal{P}_1$ of $\Psi_1$
is given for $w=(w_1,w_2) \in \mathcal{P}_2$ by the same formula as for~$\Phi$, 
namely
\begin{equation*} 
\Phi_1(w)(t) = \left\{
\begin{array}{ll}
w_1(2t),   & t \in \big[ 0,\tfrac{1}{2} \big] , \\ [0.3em]
w_2(2-2t), & t \in \big[ \tfrac{1}{2},1 \big] .
\end{array}
\right.
\end{equation*}

This construction can now be iterated as follows. Abbreviate
$(M_0,\omega_0) := (M,\omega)$
and for $n \in \mathbb{N} \cup \{0\}$ recursively define
$$
(M_{n+1},\omega_{n+1}) := (M_n \times M_n, \omega_n \oplus -\omega_n) .
$$
Note that $M_n=M^{2^n}$ is the $2^n$-fold product of~$M$ with itself. 
Let
$$
\Delta^0_1 := \Delta^1_1 := \Delta_1 \,\subset\, M_1
$$
be the diagonal in~$M_1$ and define for $n \geqslant 1$ Lagrangian submanifolds $\Delta^0_{n+1}$, $\Delta^1_{n+1}$
in~$M_{n+1}$, where the former is defined recursively as
$$
\Delta_{n+1}^0 = \Delta_n^0 \times \Delta_n^1 \,\subset\, M_{n+1},
$$
while the latter 
$$
\Delta_{n+1}^1 = \big\{(z,z) : z \in M_n \big\} \,\subset\, M_{n+1}
$$
is just the diagonal $\Delta_{n+1}$ of $M_{n+1}$.
Hence
\begin{equation} \label{e:diag}
\Delta_n^0 \,=\, \Delta_1 \times \Delta_1 \times \Delta_2 \times \Delta_3 \times \dots \times \Delta_{n-1}
\quad \mbox{for } n \geq 2.
\end{equation}

Recall that two submanifolds $S_0, S_1$ of a manifold~$X$ are said to intersect cleanly 
if $S_0 \cap S_1$ is a submanifold of~$X$ and if for every $p \in S_0 \cap S_1$
it holds that  
$$
T_p (S_0 \cap S_1) \,=\, T_p S_0 \cap T_p S_1 .
$$

\begin{lemma} \label{le:tot}
{\rm (i)}
For $n \geqslant 1$ the intersection of the two Lagrangians $\Delta_n^0$ and~$\Delta_n^1$
equals the total diagonal in~$M_{n}=M^{2^n}$, i.e.,
$$
\Delta_n^0 \cap \Delta_n^1 \,=\, \big\{(z, \ldots, z) : z \in M\big\} \,=:\, \Delta^{\tot}_n.
$$
In particular, the intersection $\Delta_n^0 \cap \Delta_n^1$ is canonically 
diffeomorphic to~$M$.

\s
{\rm (ii)}
For $n \geqslant 1$ the submanifolds $\Delta_n^0$ and~$\Delta_n^1$ intersect cleanly along $\Delta^{\tot}_n$.
\end{lemma}

\proof
(i)
We argue by induction on~$n$ starting with~$n=1$. 
For $n=1$ the assertion holds since $\Delta_1^0 \cap \Delta_1^1 = \Delta_1 \cap \Delta_1 = \Delta_1$ is the diagonal in~$M_1$. 
The inclusion $\Delta_{n+1}^0 \cap \Delta_{n+1}^1 \supset \Delta_{n+1}^{\tot}$ 
follows by induction and the definitions.
For the reverse inclusion we write an element 
$z \in \Delta_{n+1}^0 \cap \Delta_{n+1}^1$ as
$$
z=(z_1,z_2), \quad z_1,z_2 \in M_n.
$$
Since $z \in \Delta_{n+1}^1$, which is just the diagonal of~$M_{n+1}$, 
we have
$$
z_1=z_2.
$$
Since $z \in \Delta_{n+1}^0=\Delta_n^0 \times \Delta_n^1$ we have
$$
z_1 \in \Delta_n^0, \quad z_2 \in \Delta_n^1.
$$
Combining these facts we derive
$$z_1=z_2 \in \Delta_n^0 \cap \Delta_n^1$$
which by the induction hypothesis is $\Delta_n^{\tot}$. 
Using once more $z_1=z_2$ it follows that $z \in \Delta_{n+1}^{\tot}$. 

\s
(ii)
For $n=1$ the assertion is clear since $\Delta_1^0 = \Delta_1^1 = \Delta_1^{\tot}$,
and the inclusions
$$
T_z (\Delta_n^0 \cap \Delta_n^1) \,\subset\, T_z \Delta_n^0 \cap T_z \Delta_n^1
$$
for $z \in \Delta_n^0 \cap \Delta_n^1 = \Delta_n^{\tot}$ and $n \geq 1$ are also obvious.
For the reverse inclusions fix $z \in \Delta_{n+1}^{\tot}$ 
and $\zeta \in T_z \Delta_{n+1}^0 \cap T_z \Delta_{n+1}^1$
and write $z= (z_1,z_1) \in M_n \times M_n$ and $\zeta = (\zeta_1,\zeta_2) \in T_{z_1}M_n \oplus T_{z_1} M_n$.
Since $(\zeta_1,\zeta_2) \in T_z \Delta_{n+1}^1$ and $\Delta_{n+1}^1$ is the diagonal in $M_{n+1}$,
we have
$$
\zeta_1 = \zeta_2.
$$
Since $(\zeta_1,\zeta_2) \in T_z \Delta_{n+1}^0$ and $\Delta_{n+1}^0 = \Delta_n^0 \times \Delta_n^1$,
we have
$$
\zeta_1 \in T_{z_1} \Delta_n^0, \quad \zeta_2 \in T_{z_1} \Delta_n^1 .
$$
Together with the induction hypothesis it follows that
$$
\zeta_1 = \zeta_2 \,\in\, T_{z_1} \Delta_n^0 \cap T_{z_1} \Delta_n^1 \,=\, T_{z_1} \Delta_n^{\tot}.
$$
This implies that $\zeta \in T_z \Delta_{n+1}^{\tot}$.
\proofend

For $n \geqslant 1$ we define the path space
$$
\mathcal{P}_n = \big\{ v \in W^{1,2}([0,1],M_n) : 
         v(0) \in \Delta_n^0,\,\,v(1) \in \Delta_n^1\big\} .
$$
For every $n \in \NN \cup \{0\}$ the same formula as for~$\Psi = \Psi_0$ defines a 
diffeomorphism between Hilbert manifolds
$$
\Psi_n \colon \mathcal{P}_n \to \mathcal{P}_{n+1}.
$$
In particular, the following lemma holds.

\begin{lemma}
The Hilbert manifolds $\mathcal{P}_n$ for $n \in \mathbb{N} \cup \{0\}$ are all 
diffeomorphic to each other. 
\end{lemma}

The following topological lemma will imply that the solution spaces of Floer's equation
relevant for the proof of Theorem~\ref{t:main} are $C^\infty_{\loc}$-compact.
 
\begin{lemma} \label{le:omegarel}
Assume that $[\omega] |_{\pi_2 (M)} =0$.
Then $[\omega_n] |_{\pi_2(M_n,\Delta_n^i)}=0$  for $i=0,1$ and for every $n \geq 1$.
\end{lemma}

\proof
Since $\pi_2 (M_n) = \oplus_{2^n} \pi_2(M)$, we have $[\omega_n] |_{\pi_2(M_n)} =0$ for all $n \geq 1$.
The lemma now follows from the following two claims and by induction.

\m
\ni
{\bf Claim 1.}
{\it Let $(X,\omega)$ be a symplectic manifold with $[\omega] |_{\pi_2(X)} =0$, 
and let $\Delta$ be the diagonal in~$X \times X$.
Then $[\omega \oplus -\omega] |_{\pi_2 (X \times X, \Delta)} =0$.}

\s \ni
Indeed, given $u = (u_0,u_1) \colon (\mathbb{D}, S^1) \to (X \times X, \Delta)$ 
define $v \colon S^2 = \CP^1 \to X$ by $v(z) = u_0(z)$ if $|z| \leq 1$ and
$v(z) = u_1(1 / \overline z)$ if $|z| \geq 1$.
Then $\int_{\mathbb{D}} u^* (\omega \oplus -\omega) = \int_{\mathbb{D}} u_0^* \, \omega - \int_{\mathbb{D}} u_1^* \, \omega = \int_{S^2} v^* \omega =0$.

\m \ni
{\bf Claim 2.}
{\it For $i=0,1$ let $L_i \subset (X_i,\omega_i)$ be Lagrangian submanifolds of symplectic manifolds 
such that $[\omega_i] |_{\pi_2(X_i,L_i)} =0$. Then}
$$
[\omega_0 \oplus \omega_1] |_{\pi_2 (X_0 \times X_1, L_0 \times L_1)} = 0 .
$$
\ni
Indeed, given $u \colon (\mathbb{D}, S^1) \to (X_0 \times X_1, L_0 \times L_1)$ write $u = (u_0,u_1)$ 
with $u_i \colon (\mathbb{D}, S^1) \to (X_i,L_i)$.
Then
$\int_{\mathbb{D}} u^* (\omega_0 \oplus \omega_1) = \int_{\mathbb{D}} u_0^* \, \omega_0 + \int_{\mathbb{D}} u_1^* \, \omega_1 = 0+0$.
\proofend

\section{The action functional of classical mechanics} \label{s:action}

In the following discussion of the action functional of classical mechanics we assume
that $(M,\omega)$ satisfies $[\omega]|_{\pi_2(M)} =0$.

Denote by $\mathcal{P}_0^c \subset \mathcal{P}_0$
the connected component of \emph{contractible} loops in 
the free loop space $\mathcal{P}_0$ of~$M$. 
Pick a time-dependent Hamiltonian on~$M$ which depends 
periodically on time,
$$
H \in C^\infty (M \times S^1, \RR),
$$
and for $t \in S^1$ abbreviate 
$H_t = H(\cdot, t) \in C^\infty (M,\RR)$.
Given a contractible loop $v \in \mathcal{P}_0^c$ we can find 
a filling disk~$\overline{v}$ of~$v$, namely a map
$\overline{v} \in C^\infty([0,1],\mathcal{P}_0^c)$ such that
\begin{equation} \label{fildisk}
\overline{v}(1)=v
\end{equation}
and such that $\overline{v}(0)$ is a constant loop in $\mathcal{P}_0^c$. 
If $\mathbb{D} = \{ z \in \mathbb{C} : |z| \leq 1 \}$
is the unit disk in the complex plane, we can think of a filling disk as a map
$$
\overline{v} \colon \mathbb{D} \to M
$$
by setting $\overline{v}(r e^{2 \pi i t}) = \overline{v}(r)(t)$
for $r e^{2 \pi i t} \in \mathbb{D}$.
In this interpretation \eqref{fildisk} becomes
$$
\overline{v}(e^{2 \pi i t}) = v(t), \quad t \in S^1.
$$
The action functional of classical mechanics 
$\mathcal{A}_H \colon \mathcal{P}_0^c \to \RR$ is now defined by
$$
\mathcal{A}_H(v) = -\int_{\mathbb{D}} \overline{v}^* \omega-\int_0^1 H_t(v(t)) \,dt
$$
where $\overline{v}$ is a filling disk for~$v$. 
Due to the assumption that~$[\omega] |_{\pi_2(M)} =0$,
the value $\mathcal{A}_H(v)$ depends only on~$v$ but not on the choice of the filling disk~$\overline{v}$.

\m
Recall that for every $n \geqslant 1$ there is a diffeomorphism
$$
\Psi^n := \Psi_{n-1} \circ \Psi_{n-2} \circ \cdots \circ \Psi_0 \colon \mathcal{P}_0 \to \mathcal{P}_n.
$$
Abbreviate
$$
\mathcal{P}_n^c := \Psi^n ( \mathcal{P}_0^c) \,\subset\, \mathcal{P}_n
$$
the connected component of $\mathcal{P}_n$ which is the image of the contractible loops. 
Geometrically, the elements of~$\mathcal{P}_n^c$ are the paths in~$M_n$ from $\Delta_n^0$ to $\Delta_n^1$ which are homotopic through such paths to a constant path in the total diagonal $\Delta_n^{\tot} = \Delta_n^0 \cap \Delta_n^1$. 

\begin{figure}[h]
 \begin{center}
  \psfrag{0}{$0$}
  \psfrag{1}{$1$}
  \psfrag{-1}{$-1$}
  \psfrag{w}{$w$}
  \psfrag{wb}{$\overline{w}$}
  \psfrag{D0}{$\Delta_n^1$}
  \psfrag{D1}{$\Delta_n^0$}
  \leavevmode\epsfbox{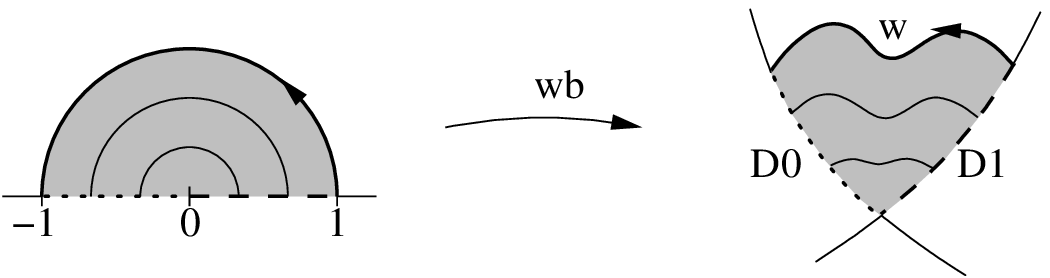}
 \end{center}
\caption{A filling half disk $\overline w$}
 \label{figure.halfw}
\end{figure}

We next discuss the pushforward of the action functional~$\mathcal{A}_H$ 
under the diffeomorphism~$\Psi^n$, i.e.,
$$
\Psi^n_* \mathcal{A}_H = \mathcal{A}_H \circ (\Psi^n)^{-1} \colon \mathcal{P}_n^c \to \mathbb{R} .
$$
We first consider the unperturbed case~$\mathcal{A}_0$, where the Hamiltonian is zero, 
which is just (minus) the area functional. For $v \in \mathcal{P}_0^c$ let $\overline{v}$ 
be a filling disk. The pushforward
$$
\overline{w} := \Psi^n \circ \overline{v} \,\in\, C^\infty([0,1], \mathcal{P}^c_n)
$$
can be geometrically interpreted as a \emph{filling half disk} of the path in~$M_n$ 
$$
w := \Psi^n \circ v \,\in\, \mathcal{P}_n^c 
$$
that starts at $\Delta^0_n$ and ends at~$\Delta^1_n$,
see Figure~\ref{figure.halfw}.
Indeed, if $\mathbb{D}_+ := \{z \in \mathbb{D} : \mathrm{Im}(z) \geqslant 0\}$ denotes the 
unit half disk lying in the upper half plane, we can think of~$\overline{w}$ as a map
$$
\overline{w} \colon \mathbb{D}_+ \to M_n
$$
by setting 
$\overline{w}(r e^{\pi i t}) = \overline{w}(r)(t)$ for $r e^{\pi i t} \in \mathbb{D}_+$.
Then $\overline{w}$ satisfies the boundary conditions
$$
\overline{w}(r) \in \Delta_n^0,\,\,r \in [0,1], \quad \overline{w}(-r) \in \Delta_n^1,\,\,r \in [0,1],\quad
\overline{w}(e^{\pi i t})=w(t),\,\,t \in [0,1] .
$$
In particular, it holds that
$\overline{w}(0) \in \Delta_n^0 \cap \Delta_n^1$.
If we now apply the pushforward of the area functional~$\mathcal{A}_0$ to $w \in \mathcal{P}_n^c$ we just obtain
$$
(\Psi^n_* \mathcal{A}_0)(w) = -\int_{\mathbb{D}_+} \overline{w}^* \omega_n
$$
where $\overline{w}$ is a filling half disk of~$w$, i.e., 
the pushforward of the area functional is the area functional again. 
For instance, for $n=1$ and with $\overline v (r,t) = \overline v (r e^{2\pi i t})$ and
$\overline w (r e^{\pi i t}) = \left( \overline w_1(r,t), \overline w_2(r,t) \right) = 
\left( \overline v(r, \frac t2), \overline v(r,1-\frac t2) \right)$, 
%
\begin{eqnarray*}
\int_{\mathbb{D}_+} \overline w^* (\omega \oplus -\omega) &=&
\int_{\mathbb{D}_+} \overline w_1^* \, \omega - \int_{\mathbb{D}_+} \overline w_2^* \, \omega \\
&=& 
\int_{\mathbb{D}_+} \overline v^*  \omega + \int_{\mathbb{D}_-} \overline v^*  \omega \\
&=& \int_{\mathbb{D}} \overline v^* \omega \;=\; - \mathcal{A}_0(v) \;=\; - \left( \Psi_*^1 \mathcal{A}_0 \right) (w) .
\end{eqnarray*}

%

We next discuss the pushforward of the Hamiltonian perturbation term
$$
\mathcal{H} \colon \mathcal{P}_0^c \to \mathbb{R}, \quad 
v \mapsto \int_0^1 H_t(v(t)) \, dt.
$$
For $n \geq 1$ define the time-dependent Hamiltonian
$$
H^n \in C^\infty(M_n \times [0,1], \mathbb{R})
$$
for $z=(z_1, z_2, \ldots, z_{2^n}) \in M_n=M^{2^n}$ and $t \in [0,1]$ by the formula
$$
H^n(z,t) = \frac{1}{2^n} \sum_{j=1}^{2^n} H\big(z_j, \tau_j^n(t) \big) .
$$
where $\tau_j^n(t)$ is recursively defined by
$\tau^0_1(t)=t$ and 
$$
\tau^{n+1}_j(t) \,=\,
\left\{\begin{array} {lcl}
\frac 12 \, \tau^n_j(t)         & \mbox{if} &  1 \leq j \leq 2^n,\\ [0.3em] 
1 - \frac 12\, \tau^n_{j-2^n}(t) & \mbox{if} &  2^n+1 \leq j \leq 2^{n+1} .
\end{array}\right.
$$
Hence
$$
\begin{array} {lll}
\tau^1_1(t) = \frac t2 & \tau^2_1(t) = \tfrac t4 & \tau^3_1(t) = \frac t8 \\ [0.4em] 
\tau^1_2(t) = 1-\frac t2 \phantom{aa}& \tau^2_2(t) = \frac 12 - \frac t4 & \tau^3_2(t) = \frac 14 - \frac t8 \\ [0.4em] 
& \tau^2_3(t) = 1 - \frac t4 & \tau^3_3(t) = \frac 12 - \frac t8 \\[0.4em] 
& \tau^2_4(t) = \frac 12 + \frac t4 \phantom{aa} & \tau^3_4(t) = \frac 14 + \frac t8 \\[0.4em] 
&& \tau^3_5(t) = 1 - \frac t8 \\[0.4em] 
&& \tau^3_6(t) = \frac 34 + \frac t8 \\[0.4em] 
&& \tau^3_7(t) = \frac 12 + \frac t8 \\[0.4em] 
&& \tau^3_8(t)= \frac 34 - \frac t8 .
\end{array}
$$

Abbreviating again 
$H^n_t = H^n(\cdot,t) \in C^\infty(M_n,\mathbb{R})$ for $t \in [0,1]$,
the pushforward of~$\mathcal{H}$ at $w \in \mathcal{P}_n^c$ becomes 
$$
\Psi^n_* \mathcal{H}(w)=\int_0^1 H^n_t(w(t)) \, dt .
$$
For instance, for $n=1$ we find using~\eqref{e:Phi},
\begin{eqnarray*}
\Psi_*^1 \mathcal{H} (w) &=& \mathcal{H} (\Phi (w)) \\
&=& \int_0^{\frac 12} H_t (w_1(2t)) \,dt + \int_{\frac 12}^1 H_t (w_2(2-2t)) \,dt \\
&=& \frac 12 \int_0^1 H_{\frac t2} (w_1(t)) \,dt + \frac 12 \int_0^1 H_{1-\frac t2} (w_2(t)) \,dt \\
&=& \int_0^1 H_t^1 (w(t)) \,dt .
\end{eqnarray*}
In particular, the pushforward of the Hamiltonian perturbation term is again a Hamiltonian perturbation term,
however of a rather specific form. Summarising we have shown that the pushforward of the action functional of classical mechanics at a point $w \in \mathcal{P}_n^c$ reads
\begin{equation}\label{pushfor}
\Psi^n_* \mathcal{A}_H(w)=-\int_{\mathbb{D}_+} \overline{w}^* \omega_n -
\int_0^1 H_t^n(w(t)) \, dt
\end{equation}
for a filling half disk $\overline{w}$ of $w$. 

Coming back to the unperturbed functionals, we have

\begin{lemma} \label{le:MB}
The area functional $\Psi_*^n \mathcal{A}_0 \colon \mcP_n^c \to \RR$ is Morse--Bott
with critical manifold the constant paths in~$\Delta_n^{\tot}$, for every $n \geq 1$.
\end{lemma}

\m \ni
{\it First proof.}
The area functional $\mathcal{A}_0 \colon \mcP_0^c \to \RR$, $v \mapsto -\int_{\mathbb{D}} \overline v^* \omega$
is Morse--Bott with critical manifold the constant loops in~$M$.
Indeed, the differential of~$\mathcal{A}_0$ at a loop $x \in \mcP_0^c$ is given by
$$
d \mathcal{A}_0(x) v \,=\, -\int_{S^1} \omega(x(t)) \left( \dot x (t), v(t) \right) dt
$$
for $v \in T_x \mcP_0^c = \left\{ v \in W^{1,2}(S^1,TM) \mid v(t) \in T_{x(t)}M \right\}$.
The non-degeneracy of $\omega$ thus implies that $\Crit \mathcal{A}_0$ is the set of constant loops in~$M$,
that we identify with~$M$.
We must show that at each $x \in M$ the kernel of the Hessian 
$d^2\mathcal{A}_0 \colon T_x \mcP_0^c \times T_x \mcP_0^c \to \RR$
is just~$T_xM$. We compute
$$
d^2\mathcal{A}_0(x) (v,w) \,=\, -\int_{S^1} \omega(x(t)) \left( \dot v(t), w(t) \right) dt .
$$
Invoking the non-degeneracy of~$\omega$ again we find that $v \in \ker d^2\mathcal{A}_0(x)$ if and only if $\dot v (t) \equiv 0$, that is, $v \in T_xM$ is constant.

Since $\Psi^n \colon \mcP_0^c \to \mcP_n^c$ is a diffeomorphism, 
the functional $\Psi_*^n \mathcal{A}_0 \colon \mcP_n^c \to \RR$ is also Morse--Bott with critical manifold the constant paths $\Psi^n (M) = \Delta_n^{\tot}$.

\m \ni
{\it Second proof.}
We fix $n \geq 1$ and abbreviate $L_i = \Delta_n^i$ for $i=0,1$.
The differential of~$\Psi_*^n \mathcal{A}_0 =: \mathcal{A}^n \colon \mcP_n^c \to \RR$, 
$v \mapsto -\int_{\mathbb{D}_+} \overline v^* \omega_n$ 
at a path $x \in \mcP_n^c$ is given by
$$
d \mathcal{A}^n(x) v \,=\, -\int_0^1 \omega_n (x(t)) \left( \dot x (t), v(t) \right) dt
$$
for $v$ in the tangent space $T_x \mcP_n^c$ of paths $v \in W^{1,2}([0,1],TM_n)$ with
$v(t) \in T_{x(t)}M_n$ and $v(i) \in T_{x(i)}L_i$ for $i=0,1$.
The non-degeneracy of~$\omega_n$ thus implies that $\Crit \mathcal{A}^n$ is the set of constant paths in~$\Delta_n^{\tot}$,
that we identify with~$\Delta_n^{\tot}$.
We must show that at each $x \in \Delta_n^{\tot}$ the kernel of the Hessian 
$d^2\mathcal{A}^n \colon T_x \mcP_n^c \times T_x \mcP_n^c \to \RR$
is just~$T_x \Delta_n^{\tot}$. 
For a constant path $x \in \Delta_n^{\tot}$,
$$
T_x \mcP_n^c \,=\, \left\{ v \in W^{1,2} ([0,1], T_xM_n) \colon v(i) \in T_xL_i \mbox{ for } i=0,1 \right\}.
$$
We compute
$$
d^2\mathcal{A}^n(x) (v,w) \,=\, -\int_0^1 \omega_n(x(t)) \left( \dot v(t), w(t) \right) dt .
$$
Invoking the non-degeneracy of~$\omega_n$ again we find that $v \in \ker d^2\mathcal{A}^n(x)$ if and only if 
$\dot v (t) \equiv 0$, that is, $v \in T_xL_0 \cap T_xL_1$ is constant.
Finally, since by Lemma~\ref{le:tot}~(ii) the two Lagrangians $L_0$ and~$L_1$ intersect cleanly along~$\Delta_n^{\tot}$, 
we have $T_xL_0 \cap T_xL_1 = T_x(L_0 \cap L_1) = T_x \Delta_n^{\tot}$.
\proofend


\section{Delay equations} \label{s:delay}

For $n \geqslant 1$ we consider a time-dependent Hamiltonian
$$
K \in C^\infty(M_n \times [0,1], \mathbb{R}) ,
$$
abbreviate $K_t = K(\cdot,t) \in C^\infty(M_n, \mathbb{R})$ for $t \in [0,1]$,
and define the action functional
$\mathcal{A}_K \colon \mathcal{P}_n^c \to \mathbb{R}$ by
$$
\mathcal{A}_K(w) = -\int_{\mathbb{D}_+} \overline{w}^* \omega_n - \int_0^1 K_t(w(t)) \, dt
$$
where $\overline{w}$ is a filling half disk for~$w$. 
Since $\Psi^n \colon \mcP_0^c \to \mcP_n^c$ is a diffeomorphism, 
the critical points of the pulled-back functional
$$
\mcA_K \circ \Psi^n \colon \mcP_0^c \to \RR
$$
are the loops $(\Psi^n)^{-1}w$ where $w$ is a critical point of~$\mcA_K$,
i.e., $w$ is a Hamiltonian chord from $\Delta_0^n$ to~$\Delta_1^n$.
In the special case that $K=H^n$ for a Hamiltonian 
$H \in C^\infty(M \times S^1, \mathbb{R})$, we have seen in~\eqref{pushfor} that
$$
\mathcal{A}_{H^n} \circ \Psi^n = \mathcal{A}_H .
$$
Therefore, in this case $(\Psi^n)^{-1}$ bijectively maps the time~$1$ $H^n$-Hamiltonian chords from $\Delta_0^n$ to~$\Delta_1^n$ in~$M_n$
to the $1$-periodic $H$-Hamiltonian orbits in~$M$.
In this section we shall see that for a general~$K$ the critical points of $\mathcal{A}_K \circ \Psi^n$ are $1$-periodic solutions of a certain delay equation on~$M$.

\begin{figure}[h]
 \begin{center}
  \psfrag{F}{$2 X_{F_{2t}}$}
  \psfrag{G}{$2 X_{G_{2-2t}}$}
  \psfrag{1}{$4 X_{F^1_{4t}}$} 
  \psfrag{3}{$4 X_{F^3_{2-4t}}$} 
  \psfrag{4}{$4 X_{F^4_{-2+4t}}$} 
  \psfrag{2}{$4 X_{F^2_{4-4t}}$} 
  \leavevmode\epsfbox{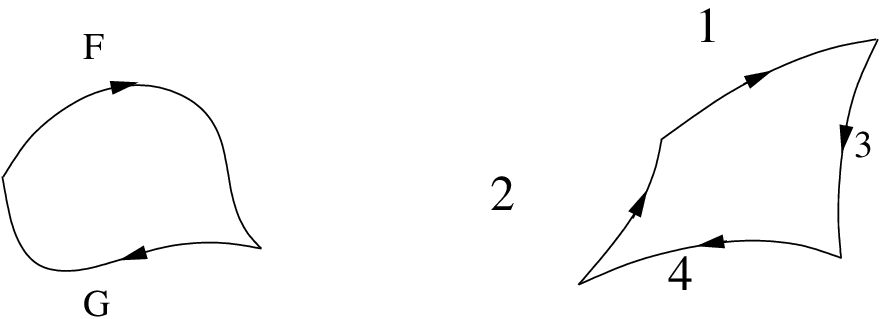}
 \end{center}
\caption{A solution of \eqref{e:sum} and of \eqref{e:sum2}}
 \label{figure.lune}
\end{figure}

\m
First take $n=1$ and recall that $M_1=M\times M$. 
Given a smooth function $K \colon M \times M \times [0,1] \to \RR$ 
define for $z_1,z_2 \in M$ the functions $M \times [0,1] \to \RR$ by
$$
K_t^{z_2}(z_1) := K_t(z_1,z_2), \quad K_t^{z_1}(z_2) := K_t(z_1,z_2).
$$
For $z = (z_1,z_2)$ we
identifying $T_z (M \times M)$ with $T_{z_1}M \times T_{z_2}M$. We can then write 
the Hamiltonian vector field~$X_K$ on $M \times M$ as
$$
X_{K_t} (z) \,=\, \left( X_{K_t}^1 (z), X_{K_t}^2(z) \right) . 
$$
Using the definition of Hamiltonian vector fields on $(M,\omega)$ and on $(M \times M, \omega \oplus -\omega)$
we obtain
\begin{equation} \label{e:X12}
X_{K_t}^1 (z) = X_{K_t^{z_2}}(z_1), \qquad X_{K_t}^2 (z) = - X_{K_t^{z_1}}(z_2). 
\end{equation}
Now assume that $w(t) = (w_1(t), w_2(t)) \in \mcP_1^c$ satisfies $\dot w (t) = X_{K_t}(w(t))$.
Let $v = \Psi^{-1}w = \Phi w \in \mcP_0^c$ be the loop in~$M$ given by formula~\eqref{e:Phi}.
Together with \eqref{e:X12} for $t \in \left[0,\frac 12 \right]$ we compute
that
\begin{eqnarray*}
\dot v(t) \,=\, \tfrac{d}{dt} (w_1(2t)) \,=\, 2 \dot w_1 (2t) &=& 2 X_{K_{2t}}^1 (w(2t)) \\
&=& 
2  X_{K_{2t}}^1 (v(t), v(1-t)) \\
&=& 2 X_{K_{2t}^{z_2=v(1-t)}} (v(t)) .
\end{eqnarray*}
In the same way  for $t \in \left[\frac 12, 1 \right]$ we find that
\begin{eqnarray*}
\dot v(t) \,=\, \tfrac{d}{dt} (w_2(2-2t)) \,=\, -2 \dot w_2 (2-2t) &=& -2 X_{K_{2-2t}}^2 (w(2-2t)) \\
&=& 
-2  X_{K_{2-2t}}^2 (v(1-t), v(t)) \\
&=& 2 X_{K_{2-2t}^{z_1=v(1-t)}} (v(t)) .
\end{eqnarray*}
Altogether, $v \in \mcP^c_0$ solves the equation
\begin{equation} \label{e:Xgen}
\dot v (t) \,=\, 
\left\{
\begin{array}{ll}
2 X_{K_{2t}^{z_2=v(1-t)}} (v(t)),    & t \in \bigl[ 0,\tfrac{1}{2} \bigr], \\ [0.3em]
2 X_{K_{2-2t}^{z_1=v(1-t)}} (v(t)) , & t \in \bigl[ \tfrac{1}{2},1 \bigr].
\end{array}\right.
\end{equation}
Note that the loop $v$ is smooth for $t \notin \{0, \frac 12\}$, but only continuous at $t = 0$ and $t=\frac 12$.
To see what kind of an equation \eqref{e:Xgen} is we specialize $K$ in several ways:

\begin{example} \label{ex:sum}
{\rm
Assume first that $K$ splits as a sum: $K_t(z_1,z_2) = F_t(z_1) + G_t(z_2)$.
Then~\eqref{e:Xgen} becomes
\begin{equation} \label{e:sum}
\dot v (t) \,=\, 
\left\{
\begin{array}{ll}
2 X_{F_{2t}} (v(t)) ,    & t \in \bigl[ 0,\tfrac{1}{2} \bigr], \\ [0.3em]
2 X_{G_{2-2t}} (v(t)) ,  & t \in \bigl[ \tfrac{1}{2},1 \bigr].
\end{array}\right.
\end{equation} 
Thus $v$ is a periodic orbit of the Hamiltonian flow defined by the
``jumping'' Hamiltonian vector field given by the right hand side of~\eqref{e:sum}: 
$v$ runs through a 2-gon, first under the law given by $2F_{2t}$, then 
under the law given by $2G_{2-2t}$, then again under the law given by $2F_{2t}$, and so on,
see the left drawing in Figure~\ref{figure.lune}.
Note that the flow defined by~\eqref{e:sum} and in particular the periodic orbit~$v$ is actually smooth
if $F_t$ and $G_t$ are constant in $z$ for $t$ near $0$ and~$1$. 
}
\end{example}

%

\begin{example} \label{ex.product}
{\rm
Assume now that $K$ splits as a product:
$K_t(z_1,z_2) = F_t(z_1)\, G_t(z_2)$.
Then 
$$
X_{K_t}^{z_2} (z_1) = G_t(z_2) \, X_{F_t}(z_1), \qquad X_{K_t}^{z_1} (z_2) = - F_t(z_1) \, X_{G_t}(z_2) ,
$$
and so $v$ solves
\begin{equation} \label{e:vFG}
\dot v (t) \,=\, 
\left\{
\begin{array}{ll}
2 G_{2t} (v(1-t)) \,X_{F_{2t}} (v(t)),       & t \in \bigl[ 0,\tfrac{1}{2} \bigr], \\ [0.3em]
2 F_{2-2t} (v(1-t)) \, X_{G_{2-2t}} (v(t)) , & t \in \bigl[ \tfrac{1}{2},1 \bigr].
\end{array}\right.
\end{equation}
This is a Hamiltonian {\it delay equation}: $\dot v (t)$ is not just given by a vector field depending
on~$t$ and $v(t)$ (as it is the case for ordinary differential equations), 
but $\dot v (t)$ also depends on the place the loop~$v$ was at another time~$1-t$.

More generally, if $K$ splits as a finite sum of products, 
$$
K_t(z_1,z_2) \,=\, \sum_{j=1}^N F^j_t(z_1)\, G^j_t(z_2) ,
$$
then $v$ solves
\begin{equation} \label{e:trigo}
\dot v (t) \,=\, 
\left\{
\begin{array}{ll}
2 \sum_{j=1}^N G^j_{2t} (v(1-t)) \,X_{F^j_{2t}} (v(t)),       & t \in \bigl[ 0,\tfrac{1}{2} \bigr], \\ [0.5em]
2 \sum_{j=1}^N F^j_{2-2t} (v(1-t)) \, X_{G^j_{2-2t}} (v(t)) , & t \in \bigl[ \tfrac{1}{2},1 \bigr].
\end{array}\right.
\end{equation}
In the case $n=1$ the general form of a delay equation we can achieve from pulling back 
Hamilton's equation for $K \colon M \times M \times [0,1] \to \RR$ by~$\Psi$ is equation~\eqref{e:Xgen}. 
Note that if $M = \RR^{2n}$ or a quotient thereof, or if $K$ is supported in a Darboux chart,  
then~\eqref{e:trigo} approximates~\eqref{e:Xgen} arbitrarily well. 
Indeed, trigonometric polynomials of period, say, $2\pi$ are linear combinations of functions of the form 
$$
F_t(z_1) \, G_t(z_2) \,=\, a(t)\bigl( e^{i (k_1,k'_1) \cdot (x_1,y_1)} \bigr) \bigl( e^{i (k_2,k_2') \cdot (x_2,y_2)} \bigr) 
$$
with $k_1,k_1', k_2,k_2' \in \ZZ$ and $a \colon [0,1] \to \RR$ a smooth function,
and trigonometric polynomials are $C^\infty$-dense in the space of all smooth periodic functions.
}
\end{example}

Recall from the introduction that classically, a delay equation on a manifold is of the form 
\begin{equation} \label{e:delayY}
\dot v(t) \,=\, \sum_{j=1}^N f_t^j \bigl( v(t-\tau_j) \bigr) \, Y_t^j  \bigl( v(t) \bigr)
\end{equation}
for time-dependent functions $f_t^j \colon M \to \RR$ and vector fields $Y_t^j$ on~$M$.
The vector fields in~\eqref{e:trigo} are not of this classical form, 
since the delay term is not taken at $v(t-\tau_j)$ but at $v(1-t)$.
As we shall see next, we can obtain classical delay equations~\eqref{e:delayY} by iterating the 
graph construction.

We start with $n=2$.
For $w = (w_1, w_2, w_3, w_4) \in \mcP^c_2$ we compute, using formula~\eqref{e:Phi} twice, that
\begin{equation} \label{e:Phi2}
v(t) = \Phi^2(w)(t) \,=\,
\left\{
\begin{array}{ll}
w_1(4t),   & t \in \big[ 0,\tfrac 14 \big], \\ [0.3em]
w_3(2-4t), & t \in \big[ \tfrac 14, \frac 12 \big], \\ [0.3em]
w_4(-2+4t), & t \in \big[ \tfrac 12, \frac 34 \big], \\ [0.3em]
w_2(4-4t), & t \in \big[ \tfrac 34, 1 \big] .
\end{array}\right.
\end{equation}
Assume now that $w$ is a Hamiltonian chord,
$\dot w = X_K(w)$ for a smooth function $K \colon M_2 \times S^1 \to \RR$.
Recall that the symplectic form on $M_2$ is $\omega \oplus -\omega \oplus -\omega \oplus \omega$.

\begin{example}
{\rm
Assume that $K$ splits as a sum,
$$
K_t(z) = F_t^1(z_1) + F_t^2(z_2) + F_t^3(z_3) + F_t^4(z_4) .
$$
Then~\eqref{e:Phi2} implies
\begin{equation} \label{e:sum2}
\dot v(t) \,=\,
\left\{
\begin{array}{ll}
4 X_{F^1_{4t}} (v(t)),    & t \in \big[ 0,\tfrac 14 \big], \\ [0.3em]
4 X_{F^3_{2-4t}} (v(t)),  & t \in \big[ \tfrac 14, \frac 12 \big], \\ [0.3em]
4 X_{F^4_{-2+4t}} (v(t)), & t \in \big[ \tfrac 12, \frac 34 \big], \\ [0.3em]
4 X_{F^2_{4-4t}} (v(t)),  & t \in \big[ \tfrac 34, 1 \big] .
\end{array}\right.
\end{equation}
Thus $v$ runs through a 4-gon, solving a different Hamiltonian differential equation
along each side,  
see the right drawing in Figure~\ref{figure.lune}.
}
\end{example}

\begin{example} \label{ex.product2}
{\rm
Assume that $K$ splits as a product,
\begin{equation} \label{e:prod2}
K_t(z) = F^1_t(z_1)\, F^2_t(z_2) \, F^3_t(z_3)\, F^4_t(z_4) .
\end{equation}
Then 
\begin{eqnarray*}
X_{K_t} (z) &=& \phantom{-} F_t^2(z_2) \, F_t^3(z_3) \, F_t^4(z_4) \, X_{F_t^1}(z_1) \\
&& - F_t^1(z_1) \, F_t^3(z_3) \, F_t^4(z_4) \, X_{F_t^2}(z_2) \\
&& - F_t^1(z_1) \, F_t^2(z_2) \, F_t^4(z_4) \, X_{F_t^3}(z_3) \\
&& + F_t^1(z_1) \, F_t^2(z_2) \, F_t^3(z_3) \, X_{F_t^4}(z_4) .
\end{eqnarray*}
To simplify notation, we set given three different numbers $j_1,j_2,j_3 \in \{1,2,3,4\}$,
$$
\mathbf{F}_t^{j_1 j_2 j_3} (v(t)) \,:=\,
F_t^{j_1} (v(\tfrac 12+t)) \cdot  F_t^{j_2} (v(\tfrac 12-t)) \cdot F_t^{j_3} (v(1-t)) . 
$$
With~\eqref{e:Phi2} and using that $v(t+1) = v(t)$ along $1$-periodic solutions,  
we then find that $v$ solves the delay equation 
\begin{equation*}
\tfrac 14 \2 \dot v (t) \,=\, 
\left\{
\begin{array}{ll}
\mathbf{F}_{4t}^{4\,32} (v(t)) \, X_{F^1_{4t}} (v(t)), & t \in \bigl[ 0,\tfrac 14 \bigr], \\ [0.5em]
\mathbf{F}_{2-4t}^{2\,14} (v(t)) \, X_{F^3_{2-4t}} (v(t)),    & t \in \bigl[ \frac 14, \tfrac 12 \bigr], \\ [0.5em]
\mathbf{F}_{-2+4t}^{1\,23} (v(t)) \, X_{F^4_{-2+4t}} (v(t)),    & t \in \bigl[ \frac 12, \tfrac 34 \bigr], \\ [0.5em]
\mathbf{F}_{4-4t}^{3\,41} (v(t)) \, X_{F^2_{4-4t}} (v(t)),    & t \in \bigl[ \frac 34, 1 \bigr] .
\end{array}\right.
\end{equation*}
%
%
%
%
%
Taking sums of Hamiltonians of the form~\eqref{e:prod2} we get more general delay vector fields.
Now take a sum of the form
$$
K_t(z) \,=\, F_t^1(z_1) \, F_t^4(z_4) + F_t^2(z_2) \, F_t^3(z_3) .
$$
Then $v$ solves the equation 
\begin{equation} \label{e:1423}
\tfrac 14 \2 \dot v (t) \,=\, 
\left\{
\begin{array}{ll}
 F^4_{4t} (v(t-\frac 12)) \, X_{F^1_{4t}} (v(t)),    & t \in \bigl[ 0,\tfrac 14 \bigr], \\ [0.5em]
 F^2_{2-4t} (v(t-\frac 12)) \, X_{F^3_{2-4t}} (v(t)),    & t \in \bigl[ \frac 14, \tfrac 12 \bigr], \\ [0.5em]
 F^1_{-2+4t} (v(t-\frac 12)) \, X_{F^4_{-2+4t}} (v(t)),    & t \in \bigl[ \frac 12, \tfrac 34 \bigr], \\ [0.5em]
 F^3_{4-4t} (v(t-\frac 12)) \, X_{F^2_{4-4t}} (v(t)),    & t \in \bigl[ \frac 34, 1 \bigr] .
\end{array}\right.
\end{equation}
The four parts of equation~\eqref{e:1423} have the form~\eqref{e:delayY} of a classical delay equation, although only with one delay, namely $\tfrac12$. }
\end{example}

For more than one delay one needs to look at higher values of $n$. We only work out the case $n=3$.

\begin{example} \label{ex.product3}
{\rm
Take $n=3$ and assume that $K \colon M_3 \times S^1 \to \RR$ splits as a product:
$$
K_t(z) \,=\, \prod_{j=1}^8 F^j_t (z_j) .
$$
Given seven different numbers $j_1, j_2, \dots, j_7 \in \{ 1,  \dots, 8 \}$ we set
\begin{eqnarray*}
\mathbf{F}_t^{j_1 \cdots j_7} (v(t)) &:=&
F_t^{j_1} (v(\tfrac 14+t)) \cdot  F_t^{j_2} (v(\tfrac 12+t)) \cdot F_t^{j_3} (v(\tfrac 34+t))\cdot  \\
&& F_t^{j_4} (v(\tfrac 14-t))  \cdot F_t^{j_5} (v(\tfrac 12-t)) \cdot  F_t^{j_6} (v(\tfrac 34-t)) \cdot F_t^{j_7}(1-t).
\end{eqnarray*}
If $w$ solves $\dot w = X_K (w)$, then $v(t) = \Phi^3(w)(t)$ solves the delay equation 
\begin{equation*}
\tfrac 18 \dot v (t) \,=\,
\left\{
\begin{array}{ll}
\mathbf{F}^{746 \, 5382}_{8t}(v(t))
\, X_{F^1_{8t}} (v(t)),    & t \in \bigl[ 0,\tfrac 18 \bigr], \\ [0.5em]
\mathbf{F}^{382 \, 1746}_{2-8t}(v(t))
\, X_{F^5_{2-8t}} (v(t)),    & t \in \bigl[ \frac 18, \tfrac 28 \bigr], \\ [0.5em]
\mathbf{F}^{461 \, 2538}_{-2+8t}(v(t))
\, X_{F^7_{-2+8t}} (v(t)),    & t \in \bigl[ \frac 28, \tfrac 38 \bigr], \\ [0.5em]
\mathbf{F}^{825 \, 6174}_{4-8t}(v(t))
\, X_{F^3_{4-8t}} (v(t)),    & t \in \bigl[ \tfrac 38, \tfrac 48 \bigr], \\ [0.5em]
\mathbf{F}^{617 \, 8253}_{-4+8t}(v(t))
\, X_{F^4_{-4+8t}} (v(t)),    & t \in \bigl[ \frac 48, \tfrac 58 \bigr], \\ [0.5em]
\mathbf{F}^{253 \, 4617}_{6-8t}(v(t))
\, X_{F^8_{6-8t}} (v(t)),    & t \in \bigl[ \frac 58, \tfrac 68 \bigr], \\ [0.5em]
\mathbf{F}^{174 \, 3825}_{-6+8t}(v(t))
\, X_{F^6_{-6+8t}} (v(t)),    & t \in \bigl[ \frac 68, \tfrac 78 \bigr], \\ [0.5em]
\mathbf{F}^{538 \, 7461}_{8-8t}(v(t))
\, X_{F^2_{8-8t}} (v(t)),    & t \in \bigl[ \tfrac 78, 1 \bigr].
\end{array}\right.
\end{equation*}
If we now take sums of products 
$$
F_t^1(z_1) F_t^4(z_4) F_t^6(z_6) F_t^7(z_7) \quad \mbox{ and } \quad 
F_t^2(z_2) F_t^3(z_3) F_t^5(z_5) F_t^8(z_8) 
$$ 
we obtain Hamiltonian delay equations with the three different time delays $\frac 14$, $\frac 12$, $\frac 34$.

More generally, pulling back Hamiltonian chords of suitable products from~$\mcP^c_n$ by~$\Psi^n$
we obtain 1-periodic solutions of Hamiltonian delay equations on~$M$ with time delays $\frac{j}{2^n}$, $j=1, \dots ,2^n-1$.
In the next section we explain how to get delay equations with different delay times.
}
\end{example}

\section{Generalization} \label{s:general}

One way to get more general delay equations on~$M$ via Hamilton's equation on iterated path spaces
is by replacing the diffeomorphism~$\Psi$, $\Psi_1, \dots$ defined in~\S \ref{s:graph}
by diffeomorphisms $\Psi \circ \sigma$, $\Psi_1 \circ \sigma_1, \dots$,
where $\sigma$ is a diffeomorphism of the circle~$S^1$ and $\sigma_1, \dots$ are 
diffeomorphisms of the interval~$[0,1]$.

Equivalently, we choose $\tau \in (0,1)$ and smooth functions $\alpha, \beta \colon [0,1] \to \RR$
such that $\dot \alpha (t) >0$ and $\dot \beta (t) <0$ and $\alpha (0)=0$, $\beta (0)=1$, $\alpha (1) = \beta (1) = \tau$,
see the figure below.
For $n \geq 0$ define $\Psi_{\alpha \beta} \colon \mcP_n \to \mcP_{n+1}$ by
\begin{equation*} 
\Psi_{\alpha \beta}(v)(t) = \left( v\big( \alpha (t) \big), v \big( \beta (t) \big) \right), \quad t \in [0,1] .
\end{equation*}
\begin{figure}[h]
 \begin{center}
  \psfrag{0}{$0$}
  \psfrag{1}{$1$}
  \psfrag{12}{$\frac 12$} 
  \psfrag{a12}{\scriptsize$\alpha_{\frac 12}$} 
  \psfrag{b12}{\scriptsize$\beta_{\frac 12}$} 
  \psfrag{a}{\scriptsize$\alpha$} 
  \psfrag{b}{\scriptsize{$\beta$}}
  \psfrag{a-}{\scriptsize$\alpha^{-1}$} 
  \psfrag{b-}{\scriptsize$\beta^{-1}$}  
  \psfrag{t}{$\tau$}  
  \leavevmode\epsfbox{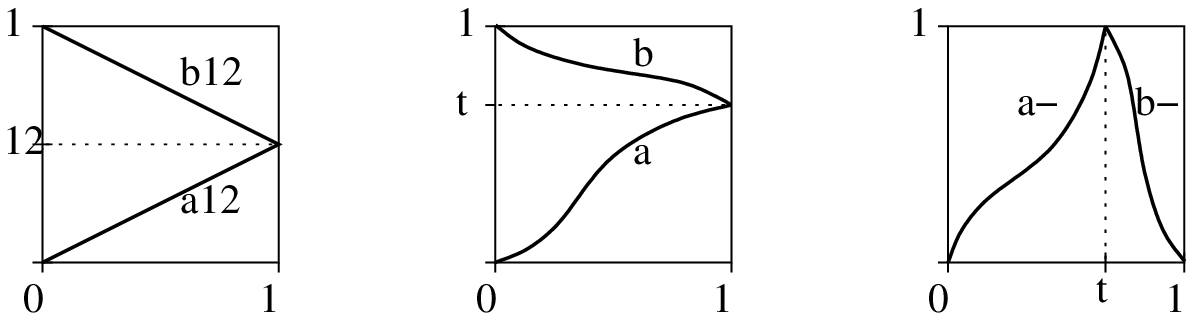}
 \end{center}
 \label{figure.ab}
\end{figure}
%

The inverse $\Phi_{\alpha \beta}$ of $\Psi_{\alpha \beta}$ is given by
\begin{equation} \label{e:Phiab}
\Phi_{\alpha \beta}(w)(t) = \left\{
\begin{array}{ll}
w_1 (\alpha^{-1}(t) ),   & t \in \big[ 0,\tau \big], \\ [0.3em]
w_2 (\beta^{-1}(t)), & t \in \big[ \tau,1 \big].
\end{array}\right.
\end{equation}

\begin{example} \label{ex:arbr}
{\rm
For $r \in (0,1)$ take 
\begin{eqnarray*}
\begin{array}{lclclcl}
\alpha_r(t) &=& rt, &&
\beta_r (t) &=& 1-(1-r)t , \\ [0.4em]
\alpha_r^{-1}(t) &=& \tfrac tr , &&
\beta_r^{-1} (t) &=& \frac{1-t}{1-r}.
\end{array}
\end{eqnarray*}
Then $\Psi_{\alpha_r\beta_r}(v)(t) = \bigl( v(rt), v(1-(1-r)t \bigr)$ and
\begin{equation*} 
\Phi_{\alpha_r \beta_r}(w)(t) = \left\{
\begin{array}{ll}
w_1 \bigl( \frac tr \bigr),   & t \in \big[ 0,r \big], \\ [0.3em]
w_2 \bigl( \frac{1-t}{1-r} \bigr), & t \in \big[ r,1 \big].
\end{array}\right.
\end{equation*}
}
\end{example}

For $n=1$ and a Hamiltonian chord $w$ of $K_t$ as in Examples~\ref{ex:sum} and~\ref{ex.product},
the loop $v = \Phi_{\alpha \beta} (w)$ solves an equation similar to~\eqref{e:sum} and~\eqref{e:vFG}.
In~\eqref{e:sum} the break time is now~$\tau$, and in~\eqref{e:vFG} the delayed times are now
$\beta \circ \alpha^{-1}(t)$ and $\alpha \circ \beta^{-1}(t)$,
which are usually different; 
for instance $\beta_r \circ \alpha_r^{-1}(t) = 1-\frac{1-r}{r}t$ and 
$\alpha_r \circ \beta_r^{-1}(t) = \frac{r}{1-r}(1-t)$ are equal only for $r=\frac 12$.
Note that in contrast to the previous delayed times, that were of the form $\pm t \pm \tau$, in these delayed times $t$ itself is scaled.
This happens in many concrete delay equations.
For instance, in the scalar non-autonomous linear pantograph equation
$$
\dot x (t) \,=\, a\, x(\lambda t) + b\, x(t)
$$
with parameter $\lambda \in (0,1)$ the delay time $\tau (t)$ is given by
$t - \tau(t) = \lambda t$, that is, $\tau (t) = (1-\lambda)t$. 

\s
Now take $n=2$. 
For $w = (w_1, w_2, w_3, w_4) \in \mcP^c_2$ we compute, using formula~\eqref{e:Phiab} twice, that
\begin{equation*} 
\Phi_{\alpha \beta}(w)(t) \,=\,
\left\{
\begin{array}{ll}
\bigl( w_1(\alpha^{-1}(t)), w_1(\alpha^{-1}(t)) \bigr) ,  & t \in \big[ 0,\tau \big] , \\ [0.3em]
\bigl( w_3(\beta^{-1}(t)), w_4(\beta^{-1}(t)) \bigr) ,    & t \in \big[ \tau, 1 \big] .
\end{array}\right.
\end{equation*}
and
\begin{equation*} 
\left( \Phi_{\alpha_1 \beta_1} \circ \Phi_{\alpha_2 \beta_2} \right) (w)(t) \,=\,
\left\{
\begin{array}{ll}
w_1 \left( \alpha_2^{-1} \circ \alpha_1^{-1}(t) \right), & t \in \big[ 0, \alpha_1(\tau_2) \big], \\ [0.3em]
w_3 \left( \beta_2^{-1} \circ \alpha_1^{-1}(t) \right), & t \in \big[ \alpha_1(\tau_2), \tau_1 \big], \\ [0.3em]
w_4 \left( \beta_2^{-1} \circ \beta_1^{-1}(t) \right), & t \in \big[ \tau_1, \beta_1(\tau_2) \big], \\ [0.3em]
w_2 \left( \alpha_2^{-1} \circ \beta_1^{-1}(t) \right), & t \in \big[ \beta_1(\tau_2), 1 \big] .
\end{array}\right.
\end{equation*}
Abbreviating
\begin{eqnarray*}
T^1(t) = (\alpha_1 \circ \alpha_2)^{-1}(t), &&
T^2(t) = (\beta_1 \circ \alpha_2)^{-1}(t), \\
T^3(t) = (\alpha_1 \circ \beta_2)^{-1}(t),  &&
T^4(t) = (\beta_1 \circ \beta_2)^{-1}(t), 
\end{eqnarray*}
we then find that if $w$ is a Hamiltonian chord,
$\dot w = X_K(w)$, for a function $K \colon M_2 \times S^1 \to \RR$ of the form
$$
K_t(z) \,=\, F_t^1(z_1) \, F_t^4(z_4) + F_t^2(z_2) \, F_t^3(z_3) ,
$$
then $v = \left( \Phi_{\alpha_1 \beta_1} \circ \Phi_{\alpha_2 \beta_2} \right) (w)$ solves the delay equation 
\begin{equation} \label{e:1423gen}
\dot v (t) \,=\, 
\left\{
\begin{array}{ll}
\phantom{-} \dot T^1(t) \,
 F^4_{T^1(t)} \bigl( v \bigl((T^4)^{-1} \circ T^1(t) \bigr) \bigr) \, X_{F^1_{T^1(t)}} (v(t)),    & t \in \bigl[ 0,\alpha_1(\tau_2) \bigr], \\ [0.5em]
- \dot T^3(t) \,
F^2_{T^3(t)} \bigl( v \bigl((T^2)^{-1} \circ T^3(t) \bigr) \bigr) \, X_{F^3_{T^3(t)}} (v(t)),    & t \in \bigl[ \alpha_1(\tau_2), \tau_1 \bigr], \\ [0.5em]
\phantom{-} \dot T^4(t) \,
F^1_{T^4(t)} \bigl( v \bigl((T^1)^{-1} \circ T^4(t) \bigr) \bigr) \, X_{F^4_{T^4(t)}} (v(t)),    & t \in \bigl[ \tau_1, \beta_1 (\tau_2) \bigr], \\ [0.5em]
- \dot T^2(t) \,
F^3_{T^2(t)} \bigl( v \bigl((T^3)^{-1} \circ T^2(t) \bigr) \bigr) \, X_{F^2_{T^2(t)}} (v(t)),    & t \in \bigl[ \beta_1 (\tau_2), 1 \bigr] .
\end{array}\right.
\end{equation}
For instance, with $\alpha_r$ and $\beta_r$ as in~Example~\ref{ex:arbr} and $T^1_{r_1r_2} = (\alpha_{r_1} \circ \alpha_{r_2})^{-1}(t)$, etc.,
the four delayed times are
\begin{eqnarray*}
\begin{array}{lclcl}
(T^4_{r_1r_2})^{-1} \circ T^1_{r_1r_2} (t) &=&
\beta_{r_1} \circ \beta_{r_2} \circ \alpha^{-1}_{r_2} \circ \alpha^{-1}_{r_1} (t) &=& \tfrac{(1-r_1)(1-r_2)}{r_1 r_2} t + r_1 \\[0.6em]
(T^2_{r_1r_2})^{-1} \circ T^3_{r_1r_2}(t) &=&
\beta_{r_1} \circ \alpha_{r_2} \circ \beta^{-1}_{r_2} \circ \alpha^{-1}_{r_1} (t) &=& \tfrac{r_2}{r_1} \tfrac{1-r_1}{1-r_2} (t-r_1) +1\\[0.6em]
(T^1_{r_1r_2})^{-1} \circ T^4_{r_1r_2} (t) &=&
\alpha_{r_1} \circ \alpha_{r_2} \circ \beta^{-1}_{r_2} \circ \beta^{-1}_{r_1} (t) &=& \tfrac{r_1 r_2}{(1-r_1)(1-r_2)} (t-r_1) \\[0.6em]
(T^3_{r_1r_2} )^{-1} \circ T^2_{r_1r_2} (t) &=&
\alpha_{r_1} \circ \beta_{r_2} \circ \alpha^{-1}_{r_2} \circ \beta^{-1}_{r_1} (t) &=& \tfrac{r_1}{r_2} \tfrac{1-r_2}{1-r_1} (t-1) +r_1 
\end{array}
\end{eqnarray*}
In particular, if $r_1=r_2=r$, then with $T^1_{r} = T^1_{rr}$, etc., 
\begin{eqnarray*}
(T^4_r)^{-1} \circ T^1_r (t) &=& \left(\tfrac{1-r}{r}\right)^2 t + r \\[0.2em]
(T^2_r)^{-1} \circ T^3_r (t) &=& t+ 1-r \\[0.2em]
(T^1_r)^{-1} \circ T^4_r (t) &=& \bigl( \tfrac{r}{1-r} \bigr)^2 (t-r) \\[0.2em]
(T^3_r)^{-1} \circ T^2_r (t) &=& t+r-1
\end{eqnarray*}
For convenience we now assume that the functions $F^j$ do not depend on time.
Taking $F^2=F^3=0$ and $r_1+r_2 = 1$, i.e., $r_1r_2 = (1-r_1)(1-r_2)$,
we then find that $v$ solves the equation 
\begin{equation} \label{e:r1r2gen}
r_1r_2 \, \dot v (t) \,=\, 
\left\{
\begin{array}{ll}
F^4 \bigl( v (t+r_1) \bigr) \, X_{F^1} (v(t)),   & t \in \bigl[ 0,r_1r_2 \bigr], \\ [0.5em]
0,    & t \in \bigl[ r_1r_2, r_1 \bigr], \\ [0.5em]
F^1 \bigl( v (t-r_1) \bigr) \, X_{F^4} (v(t)),   & t \in \bigl[ r_1, 1-r_2+r_1r_2 \bigr], \\ [0.5em]
0,    & t \in \bigl[ 1-r_2+r_1r_2 , 1 \bigr] ;
\end{array}\right.
\end{equation}
and taking $F^1=F^4=0$ and $r_1=r_2=r$ we find that $v$ solves
\begin{equation} \label{e:rrgen}
r(1-r)\, \dot v (t) \,=\, 
\left\{
\begin{array}{ll}
0,    & t \in \bigl[ 0,r^2 \bigr], \\ [0.5em]
F^2 \bigl( v (t-r+1) \bigr) \, X_{F^3} (v(t)),    & t \in \bigl[ r^2, r \bigr], \\ [0.5em]
0,    & t \in \bigl[ r, 1-r+r^2 \bigr], \\ [0.5em]
F^3 \bigl( v \bigl( t-(1-r) \bigr) \bigr) \, X_{F^2} (v(t)),    & t \in \bigl[ 1-r+r^2, 1 \bigr] .
\end{array}\right.
\end{equation}
Noting again that $v(t+r_1) = v(t-(1-r_1))$ and $v \bigl( t-r+1 \bigr) = v (t-r)$
along $1$-periodic solutions, we see that the solutions $v$ of~\eqref{e:r1r2gen} and \eqref{e:rrgen}
solve on each of their four segments a delay equation of the classical form~\eqref{e:delayY}.

\section{Proof of the Arnold conjecture for some delay equations}
\label{s:proof}

Let $(M,\omega)$ be a symplectic manifold with $[\omega] |_{\pi_2 (M)} =0$.
Fix $n \geq 1$, let $\Delta_n^0$ and~$\Delta_n^1$ be the Lagrangian submanifolds of~$(M_n,\omega_n)$
defined in \S \ref{s:graph}, and as in~\S \ref{s:action} let $\mcP^c_n$ be the space of $W^{1,2}$-paths
from $\Delta_n^0$ to~$\Delta_n^1$ that are homotopic through such paths to a constant path in the total diagonal 
$\Delta_n^{\tot} = \Delta_n^0 \cap \Delta_n^1$.
Fix a function $K \colon M_n \times S^1 \to \RR$.
Since $M$ is diffeomorphic to $\Delta_n^{\tot}$,
we must show that the number of solutions of $\dot x(t) = X_K(x(t))$ that belong to~$\mcP^c_n$
is
\begin{itemize}
\item[(i)] at least $\cuplength (\Delta_n^{\tot}) +1$;

\s
\item[(ii)]
at least $\dim H_* (\Delta_n^{\tot};\ZZ_2)$ 
if $\phi_K^1 (\Delta_n^0)$ intersects~$\Delta_n^1$ transversally.
\end{itemize}

\m \ni
{\it Proof of (ii).}
By Lemma~\ref{le:omegarel},
$[\omega_n]$ vanishes on $\pi_2(M_n,\Delta_n^i)$ for $i=0,1$,
and by Lemmas~\ref{le:tot} (ii) and~\ref{le:MB}, $\Delta_n^0$ and~$\Delta_n^1$
intersect cleanly along $\Delta_n^{\tot}$, which is the only critical manifold of the area functional
$\Psi_*^n \mathcal{A}_0$.
We can thus apply Pozniak's theorem~\cite[Theorem~3.4.11]{Po99}
to the ``isolating neighbourhood''~$\mcP_n^c$ of~$\Delta_n^{\tot}$:
The Floer homology $\HF_* (\Delta_n^0,\Delta_n^1;\ZZ_2)$ is well-defined
and isomorphic to $H_*(\Delta_n^{\tot};\ZZ_2)$. Since
$$
\HF_* (\Delta_n^0,\Delta_n^1;\ZZ_2) \,\cong\, \HF_* (\phi_K^1 (\Delta_n^0),\Delta_n^1;\ZZ_2)
$$
and since the chain complex of the latter group is generated by the Hamiltonian chords of~$K$
from $\Delta_n^0$ to $\Delta_n^1$, the claim follows.
We remark that for the special case $n=1$ we have $\Delta_1^0 = \Delta_1^1$, so that in this case
one can apply Floer's Lagrangian intersection result from~\cite{Flo88:Lag}
instead of Pozniak's theorem.

\m \ni
{\it Proof of (i).}
For the special case $n=1$, namely when $\Delta_1^0 = \Delta_1^1$, the claim follows 
from the cup-length estimate for Lagrangian intersections proved independently 
by Floer~\cite{Flo89:cup} and Hofer~\cite{Ho88}.

In general, the claim readily follows from the general approach to cup-length estimates given by Albers--Hein~\cite{AlHe16},
by taking a tubular neighbourhood~$U$ of~$\Delta_n^{\tot}$ in~$\Delta_n^0$,
and evaluating ``capped Floer strips'' with boundary on~$\Delta_n^0$ and~$\Delta_n^1$ 
at the boundary points~$(jr,0) \in \Delta_n^0$ for $j = 1, \dots, \cuplength (M)$.

\begin{figure}[h]
 \begin{center}
  \psfrag{b}{$\beta_r(s)$}
  \psfrag{s}{$s$}
  \psfrag{-1}{$-1$} 
  \psfrag{k}{$(k+1)r$} 
  \psfrag{k1}{$(k+1)r+1$}  
  \leavevmode\epsfbox{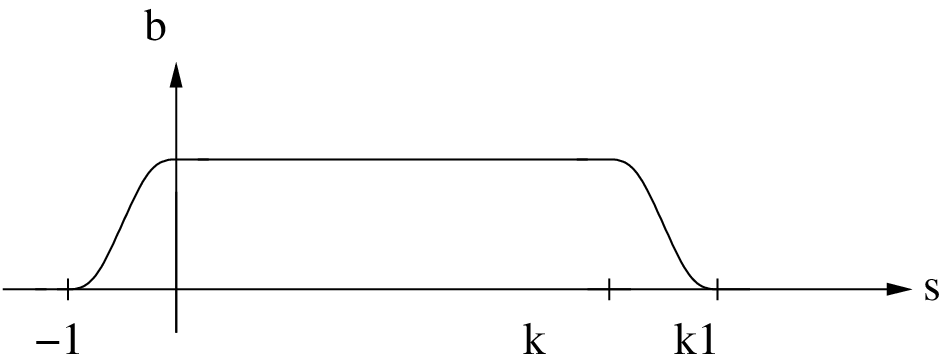}
 \end{center}
 \caption{}
 \label{figure.beta}
\end{figure}

We outline the set-up from~\cite{AlHe16}.
For this we assume that the reader is familiar with Lagrangian Floer homology.
On the path space $\mcP_n^c$ we consider the two functionals
$$
{\mathbb{A}_0}(w) = -\int_{\mathbb{D}_+} \overline{w}^* \omega_n
\quad \mbox{ and } \quad
{\mathbb{A}_K}(w) = -\int_{\mathbb{D}_+} \overline{w}^* \omega_n - \int_0^1 K_t(w(t)) \, dt .
$$
Abbreviate $k = \cuplength (M)$, and let $\beta_r$, $r \geq 0$, be a smooth family
of compactly supported functions with $\beta_0 \equiv 0$ and $\beta_r$ as in 
Figure~\ref{figure.beta} for $r \geq 1$.

%

\ni
Fix an $\omega_n$-compatible almost complex structure~$J$ on~$M_n$, and for each~$r \geq 0$ consider the the Floer equation
\begin{equation} \label{e:floer}
\pp_s u + J(u) \left( \pp_t u - X_{\beta_r(s) K(u)}\right) = 0
\end{equation}
for strips $u \colon \RR \times [0,1] \to M_n$ with $u(s,i) \in \Delta_n^i$
for $s \in \RR$ and $i=0,1$. 
Denote the space of solutions of~\eqref{e:floer} with 
$\lim_{s \to \pm \infty}u(s) \in \Delta_n^{\tot}$
by $\mathcal{M}_r$, and write $\mathcal{M}_{[0,R]} = \bigcup_{r \in [0,R]}\mathcal{M}_r$.
These spaces are all $C^\infty$-compact. 
Indeed, the only obstructions to compactness are breaking and bubbling.
From breaking one would get a non-constant $J$-holomorphic strip~$u$ asymptotic for $s \to \pm \infty$ 
to the critical manifold of~$\mathbb{A}_0$.
By Lemma~\ref{le:MB}, the only critical component of $\mathbb{A}_0$ is $\Delta^{\tot}_n = \mathbb{A}_0^{-1}(0)$.
Hence the energy of~$u$ vanishes, 
whence $u$ is constant, a contradiction.
``Bubbling'' means bubbling off of a holomorphic sphere or of a holomorphic disk at the boundary 
of a sequence of strips.
Bubbling of spheres cannot occur since $[\omega_n] |_{\pi_2(M_n)} =0$,
and bubbling of disks cannot occur since it would occur at a point in $\Delta_n^0 \cup \Delta_n^1$,
but $[\omega_n] |_{\pi_2(M_n),\Delta_n^i)} =0$ for~$i=0,1$ by Lemma~\ref{le:omegarel}.

Now choose Morse functions $f_1, \dots, f_k, f_*$ on $\Delta^{\tot}_n$, 
and extend the $f_i$ to Morse functions~$\bar f_i$ on~$M_n$ that on a tubular neighbourhood 
$U$ of~$\Delta_n^{\tot}$ in~$\Delta_n^0$ are of the form $f_i + q$
where $q \colon U \to \RR$ is a positive definite quadratic form in normal direction to~$\Delta^{\tot}_n$. 
Then $\Crit f_i \subset \Crit \bar f_i$. 

From here on the argument is exactly as in Sections~4 and~5.2 of~\cite{AlHe16}:
Fix critical points $x_*^\pm$ of~$f_*$ and $x_j$ of~$f_j$ and Riemannian metrics $g_*$ and $g_j$ on~$M_n$,
and denote by $W^u(x_*^-,f_*)$ and $W^s(x_*^+,f_*)$ the unstable manifold of~$x_*^-$ and the stable manifold of~$x_*^+$  
with respect to the negative gradient flow of~$g_*$, etc. 
For $R \geq 0$ consider the closed subspace $\mathcal{M}_R (x_1, \dots, x_k, x_*^-, x_*^+)$ of~$\mathcal{M}_R$
consisting of those $u \in \mathcal{M}_R$ with 
$$
u(-\infty) \in W^u(x_*^-,f_*), \quad u(+\infty) \in W^s(x_*^+,f_*), \quad u(jR) \in W^s(x_j,\bar f_j)
$$ 
for $j = 1, \dots, k$.
%
The mod $2$ count of zero-dimensional components of the spaces $\mathcal{M}_R (x_1, \dots, x_k, x_*^-, x_*^+)$
induces, 
after identifying Morse (co-)ho\-mo\-lo\-gy with singular (co-)homology,
operations
\begin{eqnarray*}
\Theta_R \colon H^*(\Delta_n^{\tot}) \otimes \cdots \otimes H^*(\Delta_n^{\tot}) \otimes H_*(\Delta_n^{\tot}) 
&\to& H_*(\Delta_n^{\tot}) \\
a_1 \otimes \cdots \otimes a_k \otimes b &\mapsto& (a_1 \cup \dots \cup a_k) \cap b .
\end{eqnarray*}
Since $\mathcal{M}_{[0,R]}$ is compact, the operation $\Theta_R$ does not depend on~$R$,
and since $\mathcal{M}_0$ consists of the constant maps to~$\Delta_n^{\tot}$,
$\Theta_0$ is a Morse homological realization of the cup product.
By assumption, $\Theta_0$ does not vanish, and hence all the maps $\Theta_R$ are non-zero.
Using this along a sequence $R \to \infty$ one shows as in~\cite[\S 4]{AlHe16} 
that $\mathbb{A}_K$ has at least $k+1$ critical points.

\section{Improvements}

In this paper we have only used basic results on Lagrangian Floer homology in rather simple settings. 
More sophisticated versions imply stronger forms of Theorem~\ref{t:main}.
We here describe three such improvements.

For this, the following lemma will be useful. 
Recall that a pair of submanifolds $S_0, S_1$ of a manifold~$X$ is said to be relative spin
if there exists a cohomology class $\alpha \in H^2(X;\ZZ_2)$ such that $\iota_i^* \alpha = w_2(S_i)$ for $i=0,1$.
Here, $w_2(S_i)$ is the second Stiefel--Whitney class of the tangent bundle~$TS_i$, and $\iota_i \colon S_i \to X$
are the inclusions.

\begin{lemma} \label{le:SW}
The pair $\Delta_n^0, \Delta_n^1 \subset M_n$ is relatively spin, for every $n \geq 1$. 
\end{lemma}

\proof
Abbreviate the $2^n$-fold product $1 \otimes \cdots \otimes 1 \in H^0(M_n;\ZZ_2)$ by~$1_n$,
and let $w = w_2(M) \in H^2 (M;\ZZ_2)$ be the second Stiefel--Whitney class of~$TM$.

For $n=1$, $\Delta_1^0 = \Delta_1^1 = \Delta_1$ is spin in $M \times M = M_1$.
We can take 
$$
\alpha_1 = w \otimes 1 \quad \mbox{or} \quad \overline{\alpha}_1 = 1 \otimes w .
$$
Indeed, let $\iota \colon \Delta_1 \to M \times M$ be the inclusion, 
let $\pr_1 \colon M \times M \to M$ be the projection to the first factor, 
and let $f_1 \colon \Delta_1 \to M$ be the diffeomorphism $(x,x) \mapsto x$ to the first factor.
Then $\pr_1 \circ \iota = f_1$ and $f_1^* w = w_2(\Delta_1)$. Hence
$$
\iota^* \alpha_1 \,=\, \iota^* (w \otimes 1) \,=\, \iota^* (\pr_1^* w) \,=\, f_1^* w \,=\, w_2(\Delta_1) .
$$ 
Similarly, working with the second factor of $M \times M$, we find $\iota^* \overline \alpha_1 = w_2(\Delta_1)$.

For $n=2$, $\Delta_2^0 = \Delta_1 \times \Delta_1$ and $\Delta_2^1 = \Delta_2$ are relatively spin in $M_2$.
We can take 
$$
\alpha_2 = \alpha_1 \otimes 1_1 + 1_1 \otimes \overline{\alpha}_1 \quad \mbox{or} \quad 
\overline{\alpha}_2 = \overline{\alpha}_1 \otimes 1_1 + 1_1 \otimes \alpha_1 .
$$ 
Proceeding inductively we define $\alpha_{n+1}, \overline{\alpha}_{n+1} \in H^2(M_{n+1};\ZZ_2)$ by
$$
\alpha_{n+1} = \alpha_n \otimes 1_n + 1_n \otimes \overline{\alpha}_n
\quad \mbox{and} \quad 
\overline{\alpha}_{n+1} = \overline{\alpha}_n \otimes 1_n + 1_n \otimes \alpha_n .
$$
Using~\eqref{e:diag} and elementary properties of Stiefel--Whitney classes one sees that
both $\alpha_n$ and $\overline{\alpha}_n$ define relative spin structures for $\Delta_n^0, \Delta_n^1 \subset M_n$.
\proofend

Lemma~\ref{le:SW} implies that after a perturbation (making them manifolds), 
all the moduli spaces of Floer strips used in the proofs in~Section~\ref{s:proof} 
are orientable, see~\cite[\S 8]{FOOO4} or~\cite{Schm16}.
This leads to the following improvements of Theorem~\ref{t:main}.

\m \ni
{\bf 1.}\ 
{\it The lower bounds in Theorem~\ref{t:main} can be replaced by
$\cuplength (M;\FF) +1$ and by $\dim H_*(M;\FF)$,
for every field~$\FF$.}

\proof
The Morse- and Floer moduli spaces used in the proofs of~(i) and~(ii) are orientable, 
and Pozniak's theorem holds over~$\FF$, see~\cite{Schm16}.
\proofend    

\ni
{\bf 2.}\
{\it In the situation of Theorem~\ref{t:main} assume that $n=1$.  
Then the lower bound $\cuplength (M)+1$ in assertion~(i) can be replaced by $\dim M +1$.}

\proof
The (perturbed) spaces~$\mathcal{M}_r$ of solutions of~\eqref{e:floer} with boundary on~$\Delta_1$ 
are orientable manifolds, and the same is true for the space~$\mathcal{M}_\infty$ of uncapped solutions
of equation~\eqref{e:floer} with $\beta \equiv 1$.
As in~\cite{Ho88} one finds that the evaluation map
$\ev \colon \mathcal{M}_\infty \to M \times M$, $u \mapsto u(0,0)$
induces an injection $\ev^* \colon H^*(M \times M; G) \to H^*(\mathcal{M}_\infty; G)$
for every coefficient group~$G$.
Since $[\omega] |_{\pi_2(M)} =0$, the Lusternik--Schnirelmann category of~$M$ is $\dim M +1$, see~\cite{RuOp99}. 
Hence, by the proof of Theorem~A in~\cite{Ru99}, the number of points in $\phi_K (\Delta_1) \cap \Delta_1$
is at least $\dim M +1$.
\proofend

\ni
{\bf 3.}\
Giving up the standard assumption of the paper, we
now consider a closed symplectic manifold~$M$ with $[\omega] |_{\pi_2(M)} \neq 0$.
Then the action functional~$\mathbb{A}_K$ is well-defined only on the cover of~$\mathcal{P}_n^c$ 
consisting of pairs $(w,[\overline w])$ of paths~$w \in \mathcal{P}_n^c$ together with 
a homotopy class of filling half disks~$\overline w$. 
We expect that the Floer homology of~$\mathbb{A}_K$ can still be constructed  
as a Morse--Novikov theory of this cover,
with coefficients in a Novikov ring~$\Lambda$ over~$\QQ$, 
by combining \cite{Schm16} with~\cite{FOOO17}.
This uses Lemma~\ref{le:SW}, that both $\Delta_n^0$ and~$\Delta_n^1$ are fixed point sets of (the obvious) anti-symplectic involutions,
and heavy technical machinery dealing with the occurrence of bubbling of spheres and disks.
A further complication is that 
the critical manifold~$\Delta_n^{\tot} \subset M_n$ now lifts to a disjoint union of 
copies of~$\Delta_n^{\tot}$, whence the Floer strips can now interact between different connected components of the critical manifold.

For assertion~(ii) of Theorem~\ref{t:main}, 
Floer homology should still compute the homology of~$M$, 
but with coefficients in the rational Novikov ring~$\Lambda$, 
and so assertion~(ii) with lower bound $\dim H_*(M;\QQ)$ should hold for all closed symplectic manifolds.
For $n=1$ this is a theorem, see~\cite[Theorem~1.9]{FOOO17},
and for $n \geq 2$ this should follow by combining \cite{Schm16} with~\cite{FOOO17}.

For assertion~(i), however, the cuplength estimate has to be replaced by 
a rational quantum-type cuplength estimate, 
similar to the lower bounds given in \cite{Lee.Ono96, Sch98}, 
see~\cite[Theorem~1.9]{FOOO17} for the case~$n=1$.


\end{document}